\documentclass[11pt,a4paper]{article}
\usepackage{colordvi,epsfig,amsmath,amsfonts,amssymb}
\usepackage{algorithm,algorithmic}
\usepackage{subfig}
\usepackage{url}

\usepackage{tikz}
\usetikzlibrary{arrows,positioning}
\usetikzlibrary{shapes} 

\usepackage[textsize=tiny]{todonotes}
\usepackage{gensymb} 

\usepackage{array}
\newcolumntype{H}{>{\lrbox0}c<{\endlrbox}@{}}

\newcommand{\ignore}[1]{}
\newcolumntype{H}{>{\setbox0=\hbox\bgroup}c<{\egroup}@{}}

\newcommand{\BMP}[1]{\begin{minipage}{#1\linewidth}}
\newcommand{\EMP}{\end{minipage}}
\newcommand{\HGAP}[1]{\hspace{#1\linewidth}}
\renewcommand{\i}{i} 
\renewcommand{\l}{n} 
\renewcommand{\k}{p} 
\renewcommand{\j}{m} 
\newcommand{\m}{v} 
\renewcommand{\P}{P}

\newcommand{\sonsev}{af-1}
\newcommand{\soneight}{af-2}
\newcommand{\sept}{af-3}
\newcommand{\jan}{af-4}
\newcommand{\gns}{af-5-fs}
\newcommand{\gnsorig}{af-5}
\newcommand{\nville}{af-6-fs}
\newcommand{\nvilleorig}{af-6}

\newcommand{\mcvfor}{af-7b-fs}
\newcommand{\mcvsevorig}{af-7}
\newcommand{\mcvfororig}{af-7b}

\title{On the Effectiveness of Sequential Linear Programming for the Pooling Problem}
\author{
Andreas Grothey\thanks{Email: {\tt A.Grothey@ed.ac.uk}, Homepage: {\tt http://maths.ed.ac.uk/\~{}agr/ }} \and
  Ken McKinnon
}

\date{School of Mathematics \\
  The University of Edinburgh \\
  Mayfield Road, Edinburgh EH9 3FD \\
  United Kingdom \\[.7cm]
  \today
}
\begin{document}

\maketitle

%
%
\begin{abstract}
The aim of this paper is to compare the performance of a local
solution technique -- namely Sequential Linear Programming (SLP) employing random starting points -- with state-of-the-art global solvers such as Baron and more
sophisticated local solvers such as Sequential Quadratic Programming
and Interior Point for the pooling problem. These problems can have
many local optima, and we present a small example that illustrates how
this can occur.

We demonstrate that SLP -- usually deemed obsolete since the
arrival of fast reliable SQP solvers, Interior Point Methods and
sophisticated global solvers -- is still the method of choice for an
important class of pooling problem when the criterion is the quality of
the solution found within a given acceptable time budget.
On this measure SLP significantly ourperforms all other tested algorithms.

In addition we introduce a new formulation, the qq-formulation, for
the case of fixed demands, that exclusively uses proportional
variables.  We compare the performance of SLP and the global solver
Baron on the qq-formulation and other common formulations. While Baron with
the qq-formulation generates weaker bounds than with the other
formulations tested, for both SLP and Baron the qq-formulation finds the best
solutions within a given time budget. The qq-formulation can be strengthened by pq-like cuts in which case the same bounds as for the pq-formulation are obtained. However the associated time penalty due to the additional constraints results in poorer solution quality within the time budget.

\end{abstract}

%
\ignore{
\section*{Declarations}
{\em Funding}. Both authors were self-funded.\\
{\em Conflicts of interest}. The authors declare that they have no conflict of interest.\\
{\em Availability of data, material and code}. The test problems (model and data) used in this study are available at \cite{pool:ProblemRepository}.\\
}
\section*{Acknowledgements}
We would like to thank Format International and the late Roy Fawcett at the Scottish Agricultural College for making us aware of this problem and providing us with the test data instances. 
\section{Introduction}

The Pooling Problem is the problem of mixing a set of raw materials to form a specified
set of final products in such a
way that the products satisfy a set of given limits on the
concentration of certain qualities. The composition of
these qualities in the inputs is a known parameter of the model
(although it may assumed to be stochastic in some variants).  
In the standard Diet Problem the products are directly mixed straight
from the inputs, which results in a linear model.  In the Pooling
Problems the inputs can also flow through a sets of mixing bins, as
illustrated Figure 1.  The compositions of these mixing bins are
variables of the problem and this results in the mixing constraints
being non-linear (indeed bilinear). This makes the problem non-convex
and thus it can have local solutions. 

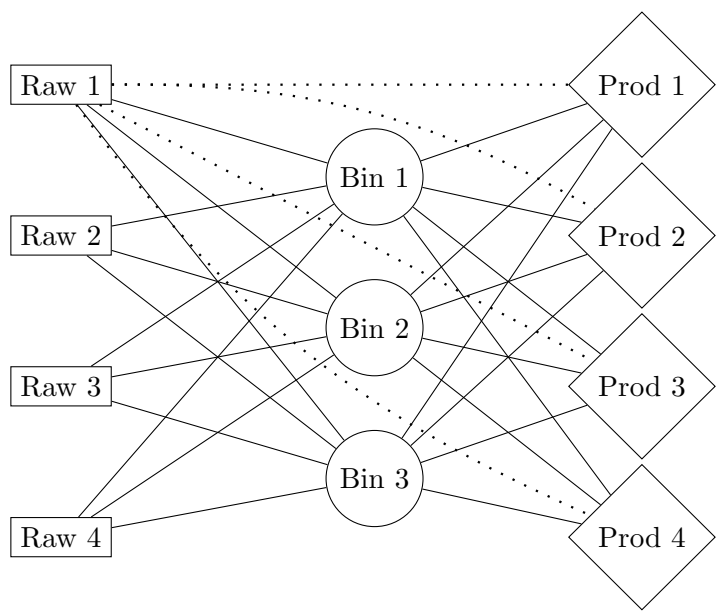
\begin{figure}[ht]
\begin{center}
\begin{tikzpicture}[node distance=2cm and 4cm,feed node/.style={rectangle,draw}, mix node/.style={circle,fill=white,draw}, product node/.style={diamond,draw}, circle dotted/.style={dash pattern=on .05mm off 8mm,line cap=round}]

\node[feed node] (F1) [] {Raw 1};
\node[feed node] (F2) [below of=F1] {Raw 2};
\node[feed node] (F3) [below of=F2] {Raw 3};
\node[feed node] (F4) [below of=F3] {Raw 4};
\node[product node] (P1) [right=6cm of F1] {Prod 1};
\node[product node] (P2) [below of=P1] {Prod 2};
\node[product node] (P3) [below of=P2] {Prod 3};
\node[product node] (P4) [below of=P3] {Prod 4};
\path[line width=0.3mm,loosely dotted] (F1) edge (P1);
\path[line width=0.3mm,loosely dotted,bend left=15] (F1) edge (P2);
\path[line width=0.3mm,loosely dotted] (F1) edge (P3);
\path[line width=0.3mm,loosely dotted, bend right=15] (F1) edge (P4);
\node[mix node] (M1) [below right=0.5cm and 3cm of F1] {Bin 1};
\node[mix node] (M2) [below of=M1] {Bin 2};
\node[mix node] (M3) [below of=M2] {Bin 3};

\path (F1) edge (M1)
           edge (M2)
           edge (M3);
\path (F2) edge (M1)
           edge (M2)
           edge (M3);
\path (F3) edge (M1)
           edge (M2)
           edge (M3);
\path (F4) edge (M1)
           edge (M2)
           edge (M3);

\path (M1) edge (P1) edge (P2) edge (P3) edge (P4);
\path (M2) edge (P1) edge (P2) edge (P3) edge (P4);
\path (M3) edge (P1) edge (P2) edge (P3) edge (P4);

\end{tikzpicture}
\caption{Superstructure of the standard Pooling Problem}
\end{center}

\end{figure}

The Pooling Problem was first described by
Haverly\cite{pool:Haverly78} in the late 1970's and shortly after
Lasdon\cite{pool:Lasdon79} proposed Sequential Linear Programming
(SLP) as a solution method.  Since then the pooling problem has become
a much studied global optimization problem with applications in the
oil and coal industry and general process optimization in chemical
engineering. It has been one of the problems driving progress in
global optimization solvers over the last few decades
\cite{pool:AndroulakisFloudas95}.

The problem of interest to us is a version of the pooling problem
arising in the modelling of animal feed mills.  
Compared to other application areas, animal feed problems are often
large scale (with up to a hundred raw materials and products and
several dozen bins) and there are many more qualities (nutrients in
this case) and restrictions on them than, say, in problems originating
from the oil or coal industry (where there are only a few qualities,
such as sulphur or octane level, that need to be taken care of). In addition the demands are fixed, firm orders, which need to be satisfied exactly, rather than maximum demand levels. Therefore, unlike the pooling formulations from the literature, the model does not have the freedom to decide which products to produce. It may seem that this substantially simplifies the problem: it is well known that the pooling problem is NP-hard\cite{pool:Haugland16}, however the proofs typically exploit this combinatorial choice to establish NP-hardness. We will show by an example that even when there are fixed orders for all products there can be many local solutions of the problem. 
Figure~1 also shows some direct connections from the raw materials to the products. These are called {\em straights}. They are typically not available for all raw material/product combinations and even where they are they can only be used at a premium cost. Further, while Figure~1 shows connections for all raw materials and bin combinations and likewise for all bin and product combinations, the set of allowable flows might be a much sparser network.
There are variants of the animal feed problem for both the standard pooling problem and the general pooling problem (with bins allowed to feed into bins), however in this paper we concentrate on the standard pooling problem.


In the context of the Pooling Problem the usual nomenclature is to
speak of inputs (or sources), pools and outputs (or targets or
products). Restrictions are on {\em qualities} of the inputs (and
outputs). We will use these terms, but in the specific context of
animal feed mills also refer to them as raw materials, (mixing) bins and
products, and we refer to the qualities as nutrients, though they can
represent more general properties such as energy or water content.

The pooling problem is a well studied global optimization problem and
most contributions to the literature are concerned with
solving the problem to proven global optimality\cite{pool:MisenerFloudas09,pool:Haugland16}. In this paper we have a
slightly different motivation: rather than solving the problem to
global optimality (which for larger problems may well not be
achievable, at least within a reasonable time) we are concerned with
finding as good a solution as possible in a limited timeframe. 
This point of view is closer to the concerns of practical
applications (at least for animal feed mills where these problems have
to be solved many times per day and additional cost due to
suboptimal solutions are not so high to make more effort worthwhile).

This paper is laid out as follows. In the following section we review the standard formulations of the pooling problem and introduce our formulation, which we term the qq-formulation, that only uses proportional variables. As far as we are aware this is a new formulation. In Section~\ref{sec:local} we provide some insight into how a small problem can have many local solutions even for fixed demands. 
Section~\ref{sec:solvers} we describe the solvers that we compare including our implementation of Sequential Linear Programming (SLP). Section~\ref{sec:results} presents numerical comparisons of the local solution methods with Baron while in Section~\ref{sec:conclusions} we draw our conclusions.

\section{Review of standard pooling problem formulations}
\label{sec:form}

The following sets, parameters and variables are used in the formulations:

\begin{itemize}
\item Sets:\\
\begin{tabular}{ll}
$\i \in \mathcal{I}$ & Set of inputs/raw materials,\\
$\j \in \mathcal{M}$ & Set of pools/bins/mixes,\\
$\k \in \mathcal{\P}$ & Set of outputs/products,\\
$\l \in \mathcal{N}$ & Set of nutrients.
\end{tabular}
\item Parameters:\\
\begin{tabular}{ll}
$r_{\i,\l}$ & nutrient composition (amount of nutrient per unit mass) of raw material $\i$,\\
$\underline{d_{\k\l}},\overline{d_{\k\l}}$ & lower and upper bounds on nutrient composition of product $\k$,\\
$c_\i$ & per unit cost of raw material $\i$ -- when used through bins,\\
$c^s_\i$ & per unit cost of raw material $\i$ -- when used directly (straights),\\
$p_\k$ & per unit selling price of product $\k$,\\
$t_\k$ & tonnages: (maximum) demand for product $\k$. 
\end{tabular}
\item Variables (p-formulation):\\
\begin{tabular}{ll}
$\m_{\j\l}$ & nutrient composition of bin $\j$,\\
$f_{\i\j}, f_{\j\k}$ & flows from raw material $\i$ to bin $\j$ and from bin $\j$ to product $\k$,\\
$f_{\i\k}$ & flow from raw material $\i$ to product $\k$ ({\em straights}).
\end{tabular}
\item Variables (other formulations):\\
\begin{tabular}{ll}
$\lambda_{\j\i}$ & proportion of bin $\j$ that originates from raw material $\i$,\\
$\mu_{\k\i}$ & proportion of demand $\k$ that originates from raw material $\i$,\\
$\mu_{\k\j}$ & proportion of demand $\k$ that originates from bin $\j$,\\
$d_{\k\l}$ & nutrient composition of demand $\k$,\\
$c^d_\k$ & per unit cost of product $\k$,\\
$c^m_\j$ & per unit cost of bin $\j$.
\end{tabular}
\end{itemize}

\subsection{PQ-formulation \textnormal{(variables $f_{\j\k}, f_{\i\k}, \lambda_{\j\i}$)}}
The earliest mathematical formulation of the problem was the {\em
  p-formulation} \cite{pool:Haverly78} which can be seen as a flow-formulation of the problem. Its variables are the total flows $f_{\i\j}$ from raw materials to bins and onto products as well as the nutrient composition of the bins $\m_{\j\l}$.

However, the standard formulation used today is the {\em pq-formulation} 
due independently to Quesada and Grossmann \cite{pool:Quesada95} and Tawarmalani and Sahinidis \cite{pool:Sahinidis02}. It is a strengthening of the earlier {\em q-formulation} proposed by Ben-Tal et al. \cite{pool:BenTal94}. Both the pq- and q-formulation 
introduce proportion variables $\lambda_{\j\i}\ge 0: \sum_\i \lambda_{\j\i} = 1$ that give the proportion of material in pool $\j$ originating from raw material $\i$, and 
expresses the nutrient content of the bins in terms of these proportional variables and the nutrient content of the raw materials. The {\em pq-formulation} of the pooling problem can be stated as
\begin{subequations}
\label{pqform}
\begin{align}
\min_{f\ge 0, \lambda\ge 0}\quad & \sum_{\i \in \mathcal{I}} \left(c_\i \sum_{\j\in\mathcal{M}}\sum_{\k\in\mathcal{\P}} \lambda_{\j\i}f_{\j\k} + c_\i^s\sum_{\k\in\mathcal{\P}} f_{\i\k}\right) - \sum_{\k \in \mathcal{\P}} p_\k\left(\sum_{\i \in\mathcal{I}} f_{\i\k} + \sum_{\j \in\mathcal{M}} f_{\j\k}\right)\label{q:obj}\\
\text{s.t.} &\notag\\
\begin{tabular}{l}product\\demand\end{tabular}\quad & \left[ \sum_{\i\in\mathcal{I}} f_{\i\k} + \sum_{\j\in\mathcal{M}} f_{\j\k} \le t_\k,\quad \forall \k\right.\label{q:proddem}\\
\text{convexity}\quad & \left[\sum_{\i\in\mathcal{I}} \lambda_{\j\i} = 1,\quad \forall \j\label{q:conv}\right.\\
\begin{tabular}{l}product\\quality\end{tabular}\quad & \left[ \sum_{\i\in\mathcal{I}} r_{\i\l}f_{\i\k} + \sum_{\j\in\mathcal{M}} \sum_{\i\in\mathcal{I}} r_{\i\l}\lambda_{\j\i}f_{\j\k}
\left\{\begin{array}{@{}l@{}l@{}}
\le & \overline{d_{\k\l}} (\sum\limits_{\i\in\mathcal{I}} f_{\i\k} + \sum\limits_{\j\in\mathcal{M}} f_{\j\k}), \\
\ge & \underline{d_{\k\l}} (\sum\limits_{\i\in\mathcal{I}} f_{\i\k} + \sum\limits_{\j\in\mathcal{M}} f_{\j\k}), 
\end{array}\right\}
\quad\forall \k,\l\right.\label{q:prodqual}\\
\text{pq-cuts}\quad &\left[\sum_{\i\in\mathcal{I}} \lambda_{\j\i} f_{\j\k} = f_{\j\k},\quad\forall \k,\j\right.
\label{q:pqcut}
\end{align}
\end{subequations}
Here (and in what follows) there are two different unit prices for each raw material: the costs of using raw-materials as {\em straights}, {\em i.e.} feeding directly into the products (through the $f_{\i\k}$ at price $c_\i^s$) is higher than when they are supplied via the mixing bins (i.e. variables $f_{\i\j}$ at cost $c_\i: c_\i < c_i^s$). Typically straights are only allowed for a subset of raw material/demand combinations.

The only bilinear terms in this formulation are the $\lambda_{\j\i}f_{\j\k}$ appearing in (\ref{q:obj}), (\ref{q:prodqual}) and (\ref{q:pqcut}): the remainder of the problem is linear.

The final set of constraints are the {\em pq-cuts}. They are redundant, indeed they are obtained by multiplying (\ref{q:conv}) with $f_{\j\k}$. Their advantage is that they provide extra (linear) constraints for the bilinear terms $\lambda_{\j\i}f_{\j\k}, \forall \i, \j, \k$ which are already part of the formulation, and this results in significantly tighter relaxations. However, their use does not come for
free: there is one of these constraints for every choice of
$\k\in\mathcal{\P}, \j\in\mathcal{M}$ which increases the number of constraints (roughly by a factor $1+|\mathcal{M}|/(|\mathcal{N}|+1)$ - see Section~\ref{sec:probsize}). The pq-formulation without the pq-cuts is indeed the q-formulation of Ben-Tal et al.

\label{pqcompvars}
The pq-formulation does not explicitly include variables giving the nutrient composition of the bins or the total use of each raw material. If there are constraints on the nutrient compositions in the bins, these can be expressed by explicitly including the variables $\m_{\j\l}$ and the linear constraints $\m_{\j\l} = \sum_{\i\in\mathcal{I}} \lambda_{\i\j}r_{\i\l}$. If constraints on the availability of raw material are given, then flows $f_{\i\j}$ from raw materials to bins would have to be calculated via the additional bilinear constraints $f_{\i\j} = \lambda_{\j\i}\sum_{\k\in\mathcal{\P}} f_{\j\k}$.

Note that in this formulation the objective
is to maximize the net profit (i.e. difference between selling price
and production cost) for each product. While there is an upper limit
$t_\k$ on how much can be produced of each product $\k$ the
optimization can decide to produce less (or even nothing at all). In
the absence of capacity restrictions (on either raw materials or bins), it is always optimal to produce
at the upper limit (if the product is to produced at all), however it
may well be optimal not to produce a product at all (in case where
this would place undue restrictions on the composition of the pools).
This adds a combinatorial choice to the other obvious non-convexities
arising from mixing constraints.

In the case of interest to us, this choice is not present: Indeed,
$t_\k$, rather than being an upper bound, is a firm order that has to be
satisfied. This has some consequences for the formulation of the
problem: constraint (\ref{q:proddem}) becomes an equality and
the terms $\sum_\i f_{\i\k} + \sum_\j f_{\j\k}$ in (\ref{q:prodqual}) and
the objective can be replaced by a constant $t_\k$. Indeed the second
half of the objective is then a constant and can be dropped.
Thus formulation (\ref{pqform}) can be replaced by the somewhat simpler form

\begin{subequations}
\label{pqform2}
\begin{align}
\min_{f\ge 0, \lambda\ge 0, d}\quad & \sum_{\i \in \mathcal{I}} \left(c_\i \sum_{\j\in\mathcal{M}}\sum_{\k\in\mathcal{\P}} \lambda_{\j\i}f_{\j\k} + c_\i^s\sum_{\k\in\mathcal{\P}} f_{\i\k}\right)\\
\text{s.t.} &\notag\\
\text{product demand}\quad & \left[ \sum_{\i\in\mathcal{I}} f_{\i\k} + \sum_{\j\in\mathcal{M}} f_{\j\k} = t_\k,\quad \forall \k\right.\label{q2:proddem}\\
\text{product quality}\quad & \left[ \sum_{\i\in\mathcal{I}} r_{\i\l}f_{\i\k} + \sum_{\j\in\mathcal{M}}\sum_{\i\in\mathcal{I}} r_{\i\l}f_{\j\k}\lambda_{\j\i} = d_{\k\l} t_\k, \quad\forall \k,\l\right.\label{q2:prodqual}\\
\text{convexity}\quad & \left[\sum_{\i\in\mathcal{I}} \lambda_{\j\i} = 1,\quad \forall \j\label{q2:conv}\right.\\
\text{pq-cuts}\quad &\left[\sum_{\i\in\mathcal{I}} \lambda_{\j\i} f_{\j\k} = f_{\j\k},\quad\forall \k,\j\right.
\label{q2:pqcut}\\
\text{bounds} & \left[ \underline{d_{\k\l}} \le d_{\k\l} \le \overline{d_{\k\l}}, \quad\forall \k, \l\right.\label{q2:bounddkl}
\end{align}
\end{subequations}
%
where the new variables $d_{\k\l}$ explicitly denote the nutrient
composition of the products. As a consequence product quality constraint (\ref{q:prodqual}) can be expressed as simple bounds (\ref{q2:bounddkl}) on the $d_{\k\l}$. Note that (\ref{q2:prodqual}) could be
used to substitute out $d_{\k\l}$ from (\ref{q2:bounddkl}) thus removing
the explicit variables. We refer to this formulation as the {\em pqs-formulation} in what follows.

\subsection{A new model: QQ-formulation \textnormal{(variables $\mu_{\j\i}, \lambda_{\k\i}, \lambda_{\k\j}, \m_{\j\l}, c^m_\j, d_{\k\l}, c^d_\k$)}}

A different formulation, which we call the {\em qq-formulation}, and
which as far as we are aware has not been described in the
literature, uses only proportion and no flow variables. That is, flows
$f_{\j\k}$ and $f_{\i\k}$ are removed from the formulation and instead
proportions $\mu_{\k\i}\ge 0, \mu_{\k\j}: \sum_\i \mu_{\k\i} + \sum_\j \mu_{\k\j} = 1$
are introduced that represent the fraction of product $\k$ that
originates from pools $\j$ or raw materials $\i$ respectively.

The nutrient composition $\m_{\j\l}$ of pools and $d_{\k\l}$ of products can be calculated using
\begin{subequations}
\begin{align}
\m_{\j\l} &= \sum_\i \lambda_{\j\i} r_{\i\l},\quad\forall \j, \l,\\
d_{\k\l} &= \sum_\i \mu_{\k\i} r_{\i\l} + \sum_\j \mu_{\k\j} \m_{\j\l},\quad \forall \k, \l.
\end{align}
\end{subequations}
Since this formulation does not include any flow variables the objective function needs to be changed. In this formulation variables $c^m_\j, c^d_\k$ representing per-unit prices of pools (mixes) and products (demands) are introduced and set via the constraints
\begin{subequations}
\begin{align}
c^m_\j &= \sum_{\i\in\mathcal{I}} \lambda_{\j\i} c_\i,\quad\forall \j\\
c^d_\k &= \sum_{\i\in\mathcal{I}} \mu_{\k\i} c_\i^s + \sum_{\j\in\mathcal{M}} \mu_{\k\j} c^m_\j,\quad \forall \k.
\end{align}
\end{subequations}
The complete qq-formulation is thus
\begin{subequations}
\label{qqform}
\begin{align}
\min_{\lambda\ge 0, \mu\ge 0, d, c^d, m, c^m}\quad & \sum_{\k \in \mathcal{\P}} t_\k c^d_\k \\
\text{s.t.} &\notag\\
\text{pool quality}\quad & \left[ \m_{\j\l} = \sum_{\i\in\mathcal{I}} \lambda_{\j\i} r_{\i\l},\quad\forall \j, \l\right.\label{qq:poolqual}\\
\text{product quality}\quad & \left[ d_{\k\l} = \sum_{\i\in\mathcal{I}} \mu_{\k\i} r_{\i\l} + \sum_{\j\in\mathcal{M}} \mu_{\k\j} \m_{\j\l},\quad \forall \k, \l\right.\label{qq:prodqual}\\
\text{convexity-m}\quad & \left[\sum_{\i\in\mathcal{I}} \lambda_{\j\i} = 1,\quad \forall \j\label{qq:conv1}\right.\\
\text{convexity-d}\quad & \left[\sum_{\i\in\mathcal{I}} \mu_{\k\i} + \sum_{\j\in\mathcal{M}} \mu_{\k\j}= 1,\quad \forall \k\label{qq:conv2}\right.\\
\text{price pools} \quad & \left[c^m_\j = \sum_{\i\in\mathcal{I}} \lambda_{\j\i} c_\i,\quad\forall \j\right.\\
\text{price products}\quad & \left[c^d_\k = \sum_{\i\in\mathcal{I}} \mu_{\k\i} c_\i^s + \sum_{\j\in\mathcal{M}} \mu_{\k\j} c^m_\j,\quad \forall \k\right.\label{qq:prodprice}\\
\text{bounds} & \left[ \underline{d_{\k\l}} \le d_{\k\l} \le \overline{d_{\k\l}}, \quad\forall \k, \l\right.
\end{align}
\end{subequations}
The bilinear terms in this formulation are $\mu_{\k\j}\m_{\j\l}$ and $\mu_{\k\j}c_\j^m$ appearing in constraints (\ref{qq:prodqual}) and (\ref{qq:prodprice}).

A constraint similar to the pq-constraint (\ref{q:pqcut}) can be
derived by substituting $\m_{\j\l}$ from (\ref{qq:poolqual}) into
(\ref{qq:prodqual}) and then multiplying (\ref{qq:conv1}) by
$\mu_{\k\j},\forall \k,\j$ to obtain
\begin{subequations}
\begin{align}
d_{\k\l} &= \sum_\i \mu_{\k\i}r_{\i\l} + \sum_\j \sum_\i \mu_{\k\j}\lambda_{\j\i}r_{\i\l},\quad\forall \k,\l\label{qq:newdkl}\\
\sum_\i \mu_{\k\j}\lambda_{\j\i} &= \mu_{\k\j},\quad\forall \k, \j\label{qq:pqcut}
\end{align}
\end{subequations}
Again the introduction of these strengthening constraints comes at the
cost of increasing the problem size. We call this strengthened qq-formulation the {\em qq+-formulation} in the later sections.

There is a close connection of the qq-formulation with the q/pq-formulations through the relations
\begin{equation}
f_{\i\k} = \mu_{\k\i}t_\k, \qquad f_{\j\k} = \mu_{\k\j}t_\k
\label{finqq}
\end{equation}
Indeed the qq-formulation can be obtained from the pq-formulation by using the above to substitute out the $f_{\i\k}$ and $f_{\j\k}$ variables. 

Note that the qq-formulation does not include any flow variables 
so if there are any capacity limits on pools or
availability limits for raw materials then flow variables would need to be added where needed by explicitly including the (bilinear) constraints (\ref{finqq}). On the other hand, limits on nutrient composition of the bins can be modelled by simple bounds on the $\m_{\j\l}$, whereas in the pq-formulation additional variables and constraints are needed to model these composition bounds as described earlier on page~\pageref{pqcompvars}. Our test problems have bounds on the nutrient composition of the bins but not on the raw material availability. The next section summarizes the situation.

\subsection{Size of Formulations\label{sec:probsize}}
The table below summarises the size (number of constraints and
variables) of the different formulations of the pooling problem as given in (\ref{pqform}), (\ref{pqform2}) and (\ref{qqform}). 
Here
N=\#nutrients, I=\#raw materials, M=\#bins, \P=\#products, and
S=\#straights (raw materials that can be used directly in products).

\begin{center}
\begin{tabular}{l|l|l}
form & variables & constraints \\
\hline
{\tt q} & $(S+M)\P + IM$ & $\P(2N\!+\!1) + M$  \\
{\tt pq} & $(S+M)\P+IM$ & $\P(2N\!+\!1) + M + M\P$ \\
{\tt pqs} & $(S+M)\P+IM +\P N$ & $\P(N\!+\!1) + M + M\P$ \\
\hline
{\tt qq} & $(S+M)\P + IM + (M+\P)(N\!\!+\!\!1)$ & $\P(N\!\!+\!\!1)+ M + M(N\!\!+\!\!1) + \P$\\
{\tt qq+} & $ (S+M)\P + IM + (M+\P)(N\!\!+\!\!1)$ & $\P(N\!\!+\!\!1) + M + M(N\!\!+\!\!1) + \P + M\P$
\end{tabular}
\end{center}
We note that
\begin{itemize}
\item The pq-constraints (\ref{q2:pqcut}/\ref{qq:pqcut}) introduce an additional $M\P$ constraints to either the q- or qq-formulation.
\item The qq-formulation has explicit $\m_{\j\l}$
  variables. When these are needed (for example to express bounds on the pool quality) the other formulations (q, pq, pqs) need a further $MN$
  variables and constraints. This accounts for
  the major size difference between the formulations. Since our test problems have bounds on the nutrient composition the size of the qq-formulations and the pqs-formulation with those added variables and constraints are roughly comparable. 
\item  If there are limits on the amount of available raw materials, these can be expressed directly as $I$ linear constraints in formulations with explicit flow variables (namely $p$, $q$ and $pq$), while the qq-formulation would need to introduce additional bilinear constraints. 
\item The qq-formulation has an additional $M+P+PN$ variables and $M+P$
  constraints compared to the q/pq/pqs-formulations. These are due to the
  explicit $c^m, c^d$ and $d_{\k\l}$ variables (the latter two of which could be
  substituted out).
\end{itemize}


\section{Occurrence of Local Solutions\label{sec:local}}

Despite having fixed demands the pooling problems that arise in animal feed mills
often have large numbers of local optima.  To understand the reason for this it
is helpful to view the problem in nutrient space, with the compositions of each
raw material, product and bin as points in this space.  The bins must lie in the
convex hull of the raw materials, and for the problem to be feasible the product
specifications must also lie there.

\begin{center}
\begin{figure}[hbt]
\includegraphics[width=0.65\textwidth]{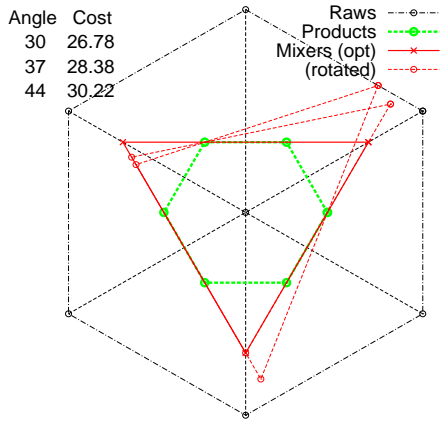}
\caption{Optimal and forced rotations}
\label{Rotate}
\end{figure}
\end{center}

Fig.~\ref{Rotate} shows an illustrative example with 2 nutrients.  The two
dimensions are the amounts of each nutrients per unit weight.  There are 7 raw
materials, 6 at the vertices of the outer black hexagon, each with a unit cost
of 6, and one at its centre, with a unit cost of 1. There are no limits on their
supply. There are 6 products, which are at the corners of the inner green
hexagon, each with demand of 1, and 3 mixer bins shown in red, whose composition
depends on the amounts of raw materials supplying it. Since the central raw
material is cheaper than the others, the unit cost of the mixture in a bin
increases with its distance from the centre.
Consider first the case where no straights are used. To be feasible the convex
hull of the bins must contain all the products. The optimization problem can
therefore be viewed as finding the bin triangle that contains all the products
that is as close to the centre as possible. The solid red triangle in
Fig.~\ref{Rotate} shown one global optimum. There are another five symmetric global
optima with the bin triangles rotated through $60\degree$, however all
intermediate positions are worse.  Two sub-optimal solutions are shown. These are
the best possible configurations where a bin is forced to lie at angles
$37\degree$ or $44\degree$. Forcing the rotation forces the average position of
the bins to move outwards, so increasing the cost.

\begin{figure}[hbt]
\HGAP{-0.12}
\BMP{0.33}
\includegraphics[width=1.58\textwidth]{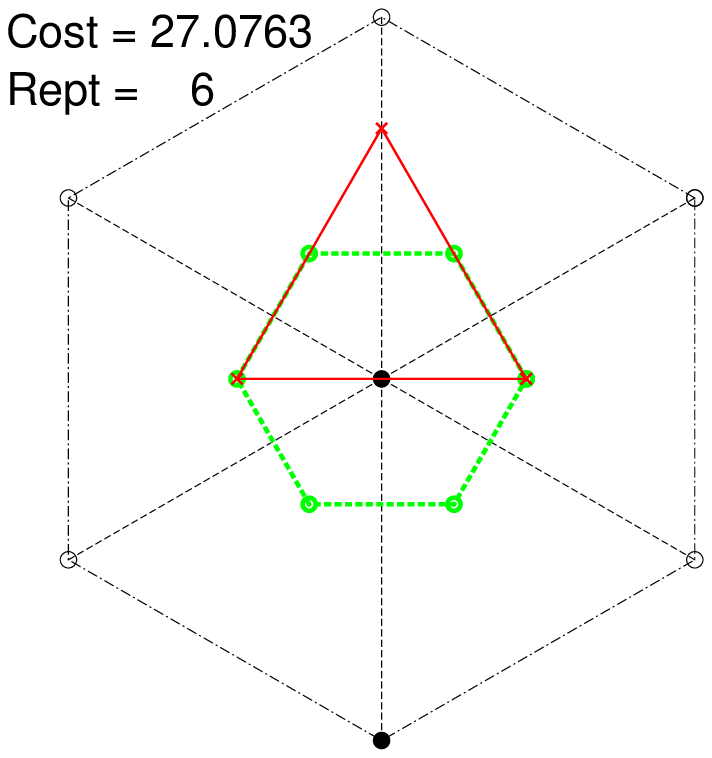}
\EMP\HGAP{0.01}
\BMP{0.33}
\includegraphics[width=1.58\textwidth]{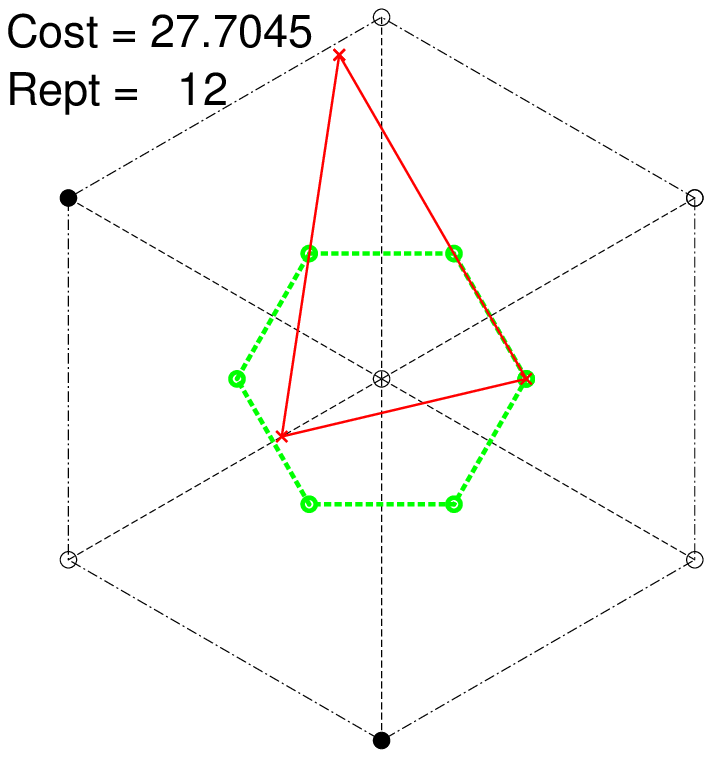}
\EMP\HGAP{0.01}
\BMP{0.33}
\includegraphics[width=1.58\textwidth]{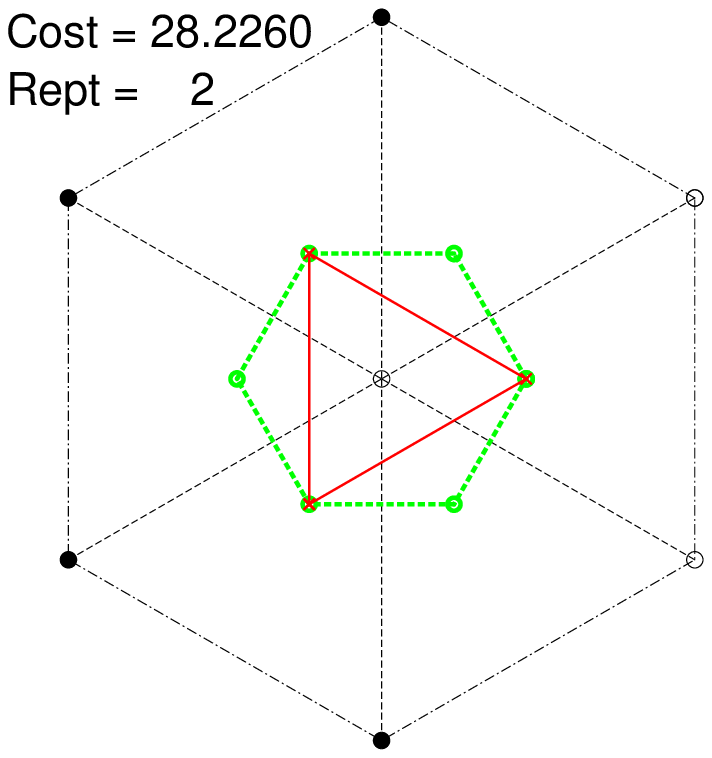}
\EMP

\HGAP{-0.12}
\BMP{0.33}
\includegraphics[width=1.58\textwidth]{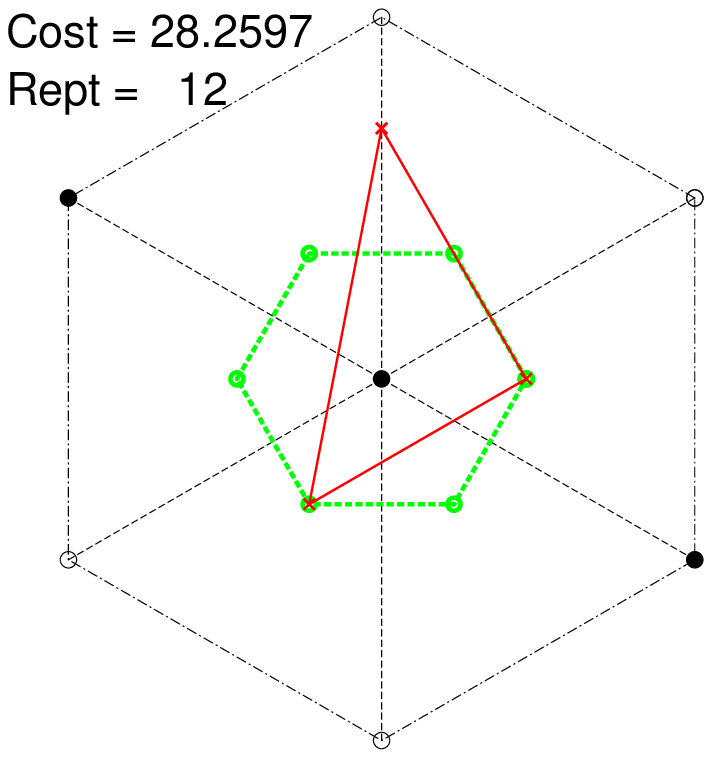}
\EMP\HGAP{0.01}
\BMP{0.33}
\includegraphics[width=1.58\textwidth]{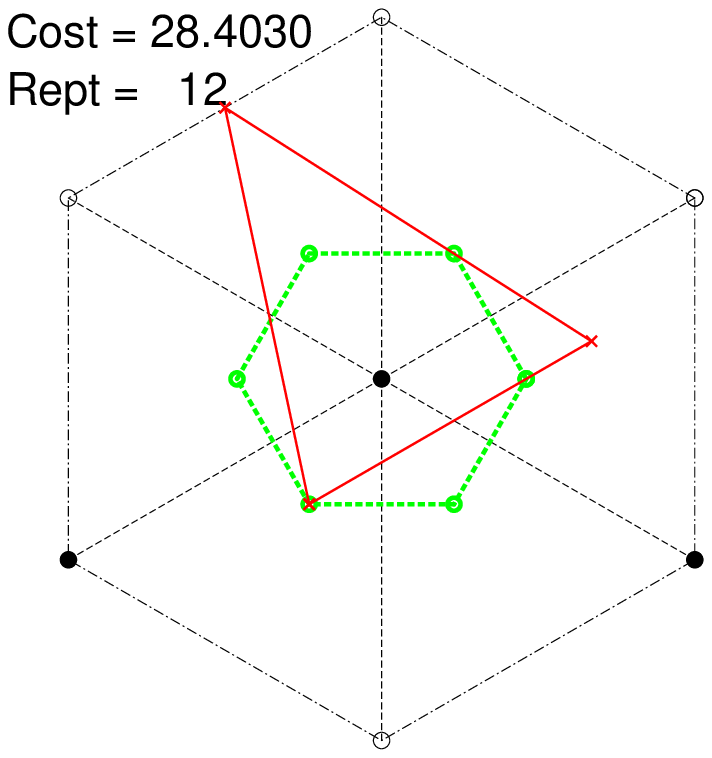}
\EMP\HGAP{0.01}
\BMP{0.33}
\includegraphics[width=1.58\textwidth]{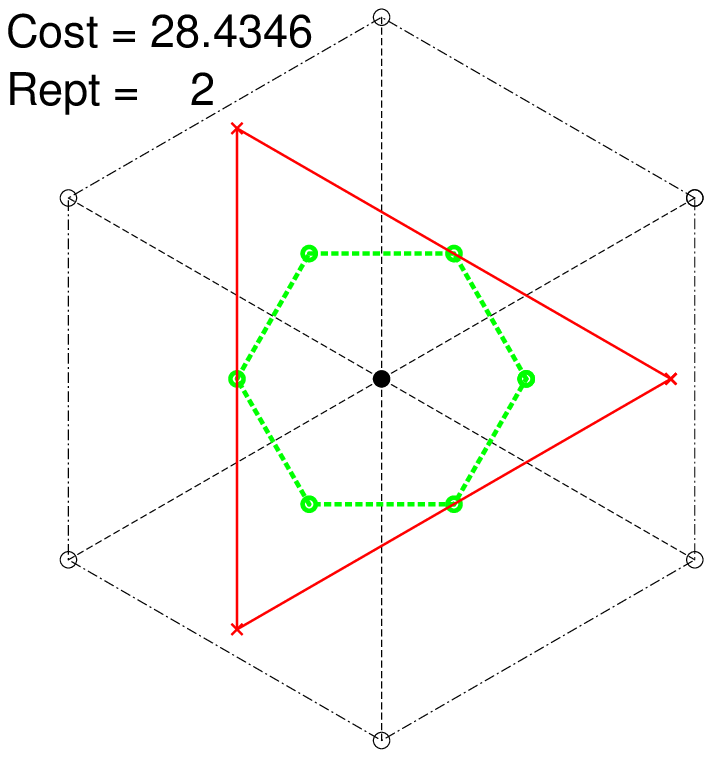}
\EMP

\HGAP{-0.12}
\BMP{0.33}
\includegraphics[width=1.58\textwidth]{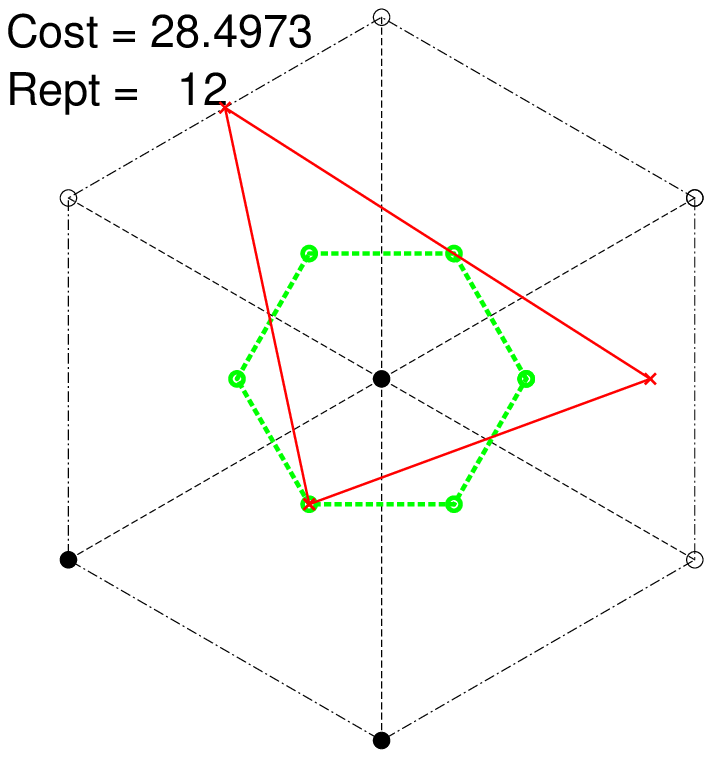}
\EMP\HGAP{0.01}
\BMP{0.33}
\includegraphics[width=1.58\textwidth]{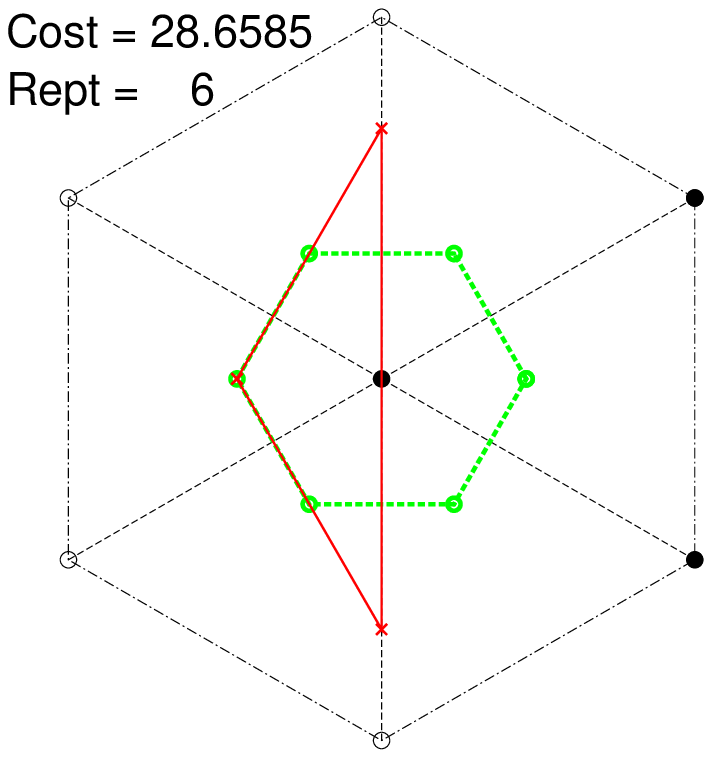}
\EMP\HGAP{0.01}
\BMP{0.33}
\includegraphics[width=1.58\textwidth]{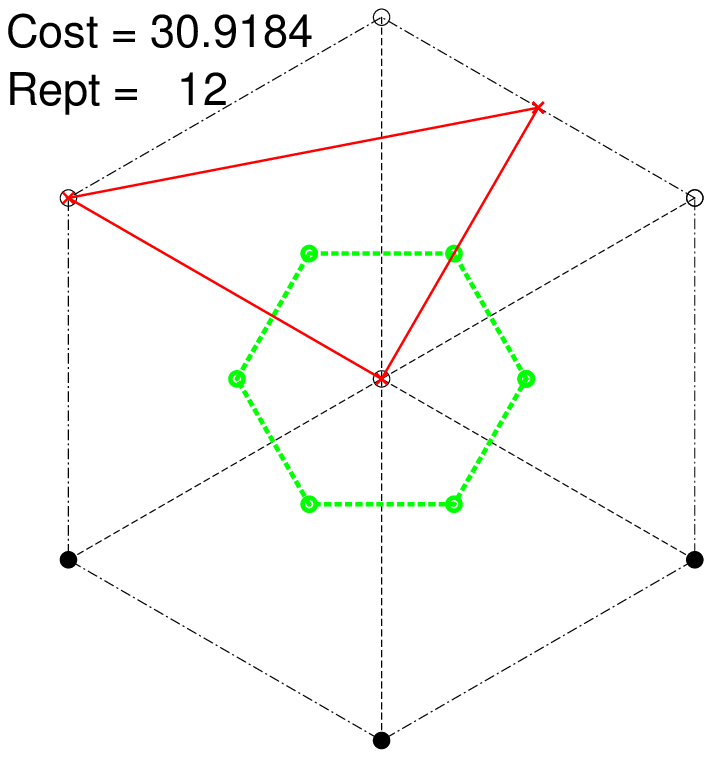}
\EMP

\HGAP{-0.12}
\BMP{0.33}
\includegraphics[width=1.58\textwidth]{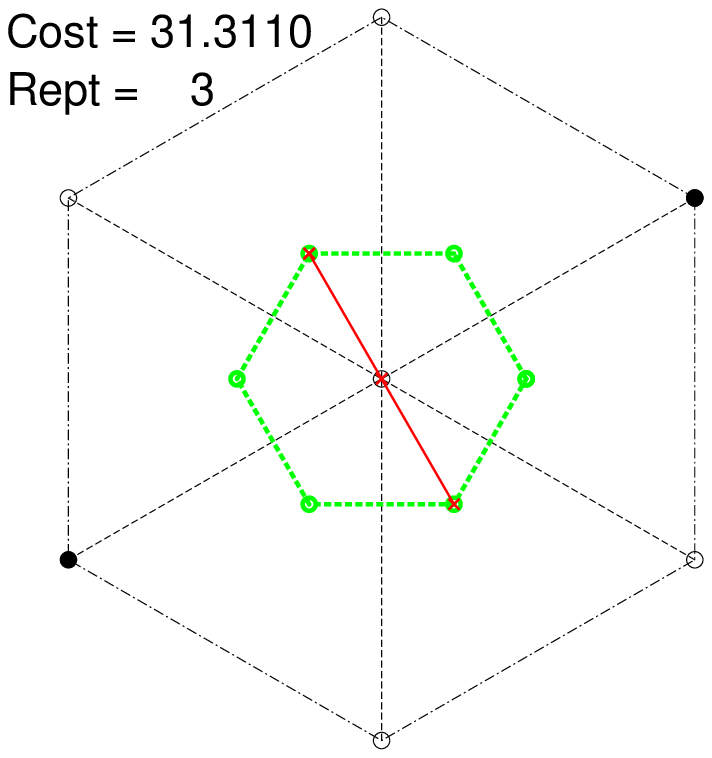}
\EMP\HGAP{0.01}
\BMP{0.33}
\includegraphics[width=1.58\textwidth]{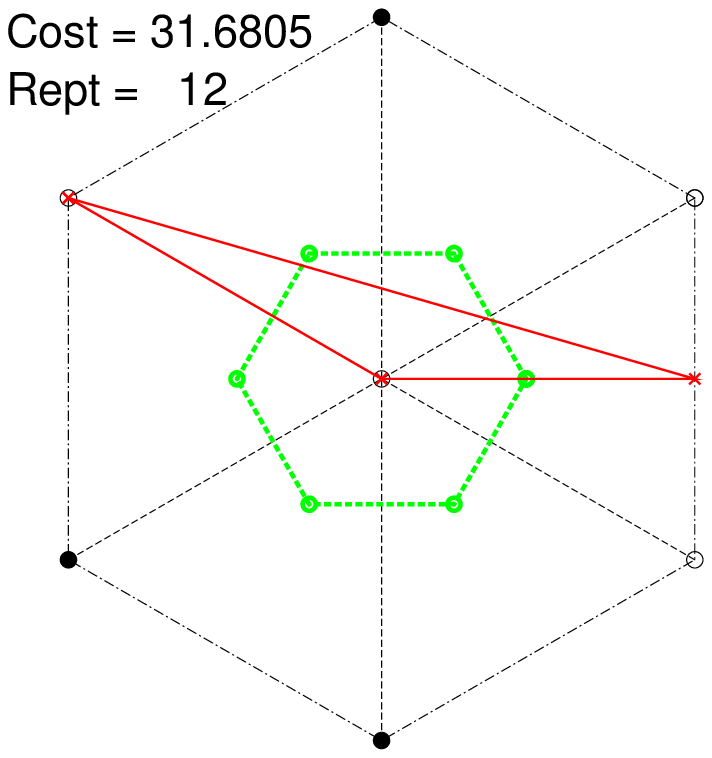}
\EMP\HGAP{0.01}
\BMP{0.33}
\includegraphics[width=1.58\textwidth]{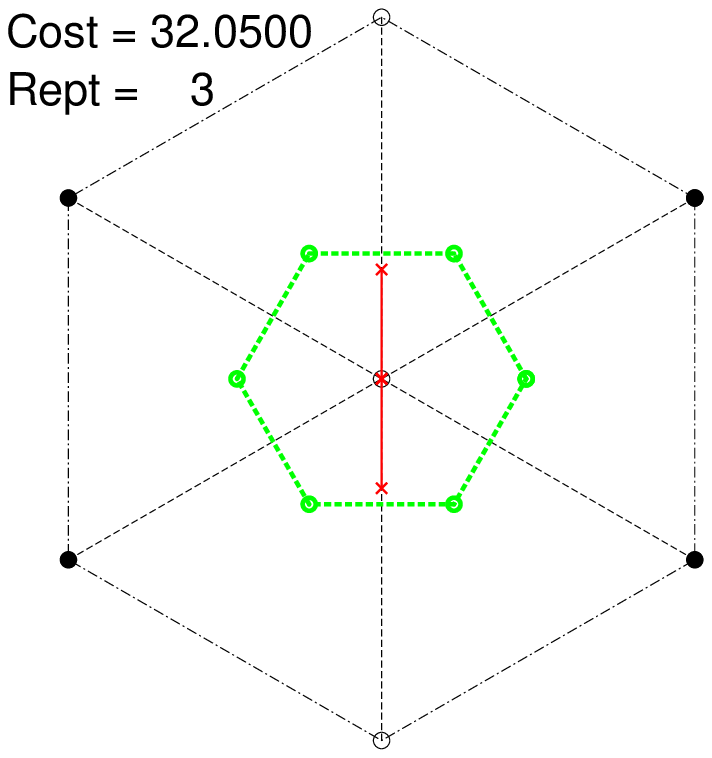}
\EMP

\caption{Local Optima: All except the final 4 cases are isolated local optima. The final 4 cases are one example from a 1- or 2-dimensional flat area of local optima.}
\label{LocalSols}
\end{figure}

Now consider the case where it is possible to supply raw materials straight to
products without passing thought the mixer bins, but at an extra cost.  This
removes the need for the convex hull of the bins to contains all the products,
and can lead to local optima where the convex hull containing different subsets
of the products. Fig~\ref{LocalSols} shows all the local solutions for the case
when the cost of straights is 2.1 time the cost of supply via the bins. We
report the cost and the number of repetitions of this solution due to
symmetries. When a raw is used as a straight its location is shown with a black
dot.  In total there are 94 local optimal in addition to the six global
solutions in Fig~\ref{Rotate}. Due to symmetry in the example there are only 13
different values of the cost, but a minor perturbation of the demands or raw
material costs would remove this symmetry without destroying the local
optimality, and in that case there would be 96 distinct objective values.  There
is an additional symmetry in the problem as we are assuming all the bins are
interchangeable.  If this is not the case, for example because the have different
capacities or are in different locations with different transport costs,
then the number of solutions could increase by a factor of 6.

Finally note that in this example the global solution was the one that did not use any
straights. This is due to the relatively high cost of straights. At lower costs
this is no longer the optimal solution, and other local solutions become the global one.


\section{Solution Methods\label{sec:solvers}}

The first solution method proposed for the pooling problem was
Sequential Linear Programming (SLP) by Lasdon\cite{pool:Lasdon79}.
SLP, as a local method, is not guaranteed to converge to a global minimizer and may even terminate at a local minimum of the corresponding feasibility problem.
Even as a local solution method for nonlinear programming problems,
SLP has been deemed obsolete due to the development of more sophisticated methods such as Sequential Quadratic Programming
(SQP,\cite{pool:FletcherLeyffer02}) and Interior Point Methods (IPM,
\cite{pool:Wright97}). Progress made in the implementations of these
methods have made the local solution of pooling problems very fast.

More recently advances in dealing with bilinear constraints in MINLP
methods \cite{pool:MisenerFloudas09} within Outer Approximation Branch
\& Bounds solvers such as Couenne and Baron, have made the global solution
pooling problem tractable at least for small instances.

The aim of this study is to evaluate the performance of different local solution methods, namely SQP, SLP and Interior Point using random multistart, and to compare this with a global solution method.

\paragraph{Sequential Linear Programming (SLP).}
In order to solve a constrained nonlinear optimization problem
\begin{equation}
\min_x f(x),\quad \text{subject to~} g(x)\le 0,
\label{nlp}
\end{equation}
the basic sequential linear programming method uses successive
linearizations of the problem around the current solution estimate
$x^{(k)}$. That is, given $x^{(k)}$, SLP
solves the problem
\begin{equation}
\min_{\Delta x} f(x^{(k)}) + \nabla f(x^{(k)})^T \Delta x, \quad
\text{s.t.~} \nabla g(x^{(k)})^T \Delta x \le -g(x^{(k)}).
\label{SLPstepv1}
\end{equation}
and then takes a step $x^{k+1} = x^k + \Delta x$. 
Such an iteration will converge to a local solution of the problem
only under fortuitous circumstances: 
even if $x^{(k)}$ is close to the optimal solution, problem (\ref{SLPstepv1}) may be unbounded or lead to very large steps.\footnote{Local convergence from starting points close to the solution is typically only ensured when the solution to (\ref{nlp}) is at a vertex of the constraints.}
In practice the SLP subproblem is therefore wrapped in a
Trust Region scheme\cite{pool:ConnGouldTointTRBook}, that is, given a trust region radius
$\rho^{(k)}>0$, we solve
\begin{equation}
\min_{\Delta x} f(x^{(k)}) + \nabla f(x^{(k)})^T \Delta x, \quad
\text{s.t.~} \nabla g(x^{(k)})^T \Delta x \le -g(x^{(k)}), \|\Delta
x\|_\infty \le \rho^{(k)}
\label{SLPstepv2}
\end{equation}
and employ the usual Trust Region methodology: that is after every
solution we compare the improvement in function value and constraint
violation predicted by the linearized model with what can actually be
achieved by taking the step $x^{(k)} + \Delta x$. Depending on the
outcome of the test we either choose to take the step, that is
$x^{(k+1)} = x^{(k)} + \Delta x$ (and possibly enlarge the trust
region), or we reject the step, reduce the trust region and keep the
same iterate $x^{(k+1)} = x^{(k)}$. Since the test of goodness of the
step $\Delta x$ consists of two criteria (objective function and
constraint violation) either a merit function or a
filter\cite{pool:FletcherLeyffer02} can be employed. In our implementation we
are using a filter strategy. 

\paragraph{Sequential Quadratic Programming (SQP)} essentially uses the same
methodology, but rather than solving the linear approximation
(\ref{SLPstepv2}) it augments this problem by adding the Hessian of
the Lagrangian to the problem: that is given a primal-dual estimate
$(x^{(k)}, \lambda^{(k)})$ of the optimal solution and the constraint
multipliers at that point, SQP solves the problem 
\begin{align}
\min_{\Delta x} & f(x^{(k)}) + \nabla f(x^{(k)})^T \Delta x +
\tfrac{1}{2}\Delta x^T \nabla^2_{xx}\mathcal{L}(x,\lambda) \Delta x, \label{SQPstepv2}\\
\text{s.t.~} & \nabla g(x^{(k)})^T \Delta x \le -g(x^{(k)}), \notag\\
& \|\Delta x\|_{\infty} \le \rho^{(k)}\notag
\end{align}
where $\mathcal{L}(x,\lambda) = f(x) - \sum_{i=1}^{m} \lambda_i g_i(x)$ is the Lagrangian function of problem (\ref{nlp}).

The second order term can be
motivated by realising that the step calculated by
(\ref{SQPstepv2}), without the trust region constraint is the same
step that would be calculated by Newton's Method employed to find a stationary point of the Lagrangian $\mathcal{L}(x, \lambda)$. See for example 
\cite{pool:FletcherLeyffer02} for further details on the implementation of an SQP method.

\section{Results\label{sec:results}}

To evaluate the efficiency of the various problem formulations and solution methods we have
tested them on a set of problems taken from the animal feed mix
industry. 

We have used compared four different formulations, namely the standard pq-formulation ({\tt pq}), the pq-formulations including the simplifications due to fixed demands ({\tt pqs}), the qq-formulation ({\tt qq}) and the qq-formulation  strengthened with the pq-equivalent cuts (\ref{qq:pqcut}) ({\tt qq+}). The solvers that we have compared are Baron as a global solver and our own implementation of SLP, FilterSQP\cite{pool:FletcherLeyffer02} and IPOpt (as an Interior Point Method).

Our test problems and their sizes are summarised in
Table~\ref{tab:testproblems}. These are problem instances from industrial practice. In most problems only a subset of the raw materials can be used as straights (i.e. directly feeding into the products). The variants with names ending on `fs', however, allow the full set of straights (at a cost of $10\times$ the normal raw material cost). The problem set is available for download from\cite{pool:ProblemRepository}.  The final two columns in Table~\ref{tab:testproblems} give
the best objective value that has been found by any method in our
tests and the best lower bound found by Baron, within 2h of computation time, for any of the tested
formulations. Note that for none of the problems and formulations was
Baron able to prove optimality. For the strengthened formulations ({\tt pq, pqs, qq+})
the gap was between 0.1\% (\mcvfor) and 4.5\% (\nvilleorig), but was between 59\% up to 100\% (lower bound of $0$)
for formulation {\tt qq}.

\begin{table}[hb]
\begin{center}
\begin{tabular}{l|rrrrr|rr|rrr}  
problem   & N & I & M & P &  S & n  & m  &  best known obj & best LB & gap\\     
\hline
\sonsev &  8 &16 & 6 &25  &16 & 925 & 310 &296.44604 & 295.93153 & 0.1\%\\
\soneight &  9 &35 & 7 &30  & 2 & 871 & 407 & 239334.02 & 235751.11 & 1.5\%\\ 
\sept   & 11 &31 & 7 &27  & 2 & 819 & 442 &122257.45 & 120498.73& 1.4\%\\
\jan    & 14 &28 & 6 &32  & 4 & 973 & 608 &186097.11 & 183620.39 & 1.3\%\\
\gns      & 17 &29 & 7 &18  &29 &1266 & 475 &119016.92 & 118763.29 & 0.2\%\\
\gnsorig& 17 &29 & 7 &18  & 5 & 816 & 475 &124003.48 & 123243.60 & 0.6\%\\
\nville   & 17 &31 & 7 &40  &31 &2555 & 893 &  2921.30 & 2881.41 & 1.4\%\\
\nvilleorig& 17 &31 & 7 &40  & 4 &1402 & 893 &  3074.23 & 2937.68 & 4.5\%\\
\mcvfor    & 14 &35 &14 &50  &35 &3788 &1024 &130290.06 & 130147.66 & 0.1\%\\
\mcvsevorig  & 14 &35 & 7 &50  & 8 &1690 & 912 &151381.81& 145928.01 & 3.6\%\\
\mcvfororig & 14 &35 &14 &50  & 8 &2334 &1024 &147285.00& 145928.01 & 0.9\%
\end{tabular}
\caption{Test problem statistics: N=\#nutrients/qualities, I=\#raw materials/inputs, M=\#mixing bins/pools, P=\#products/outputs, S=\#straights (inputs that are allowed to feed directly to outputs), n, m = \#variables and constraints in the qq-formulation.}
\label{tab:testproblems}
\end{center}
\end{table}

\subsection{Local Solvers: SLP, SQP and Interior Point}
To compare the different solution algorithms we start with the local
solvers: namely SLP, SQP and Interior Point. 
We have used the {\tt qq}-formulation for all these runs since it was observed to perform best.

As the Interior Point solvers in these comparisons we have used
IPOpt\cite{pool:WaechterBiegler06}, for SQP we have used
FilterSQP\cite{pool:FletcherLeyffer02} and our own implementation of a
filter-SLP algorithm. FilterSQP uses the active set solver {\tt
  bqpd}\cite{pool:Fletcher93} as the QP solver, 
whereas SLP uses 
CPLEX (primal Simplex, which was found to work best in this setting)
as the LP solver. Both of these employ LP/QP hotstarts between SLP/SQP iterations and a pre-solve phase that performs bound tightening and scaling of variables. Although {\tt bqpd}
uses a sparse linear algebra implementation, CPLEX is
significantly faster than bqpd when both are used to solve the same
LP. In order to give a fair comparison of SLP against SQP independent
of the subproblem solver used we have in fact tested three SLP/SQP
solvers: FilterSQP, FilterSQP as an SLP solver (by simply passing it
zero Hessians), and SLP-CPLEX.
Apart from the different LP solvers, SLP-CPLEX when compared to FilterSQP without Hessians uses a more aggressive trust region logic
and the removal of all features that make use of Hessian information
(such as second order correction steps).
As we show below, even the ad-hoc SLP setup in FilterSQP-noHess
shows some of the advantages of SLP vs SQP, whereas SLP-CPLEX is
significantly better than either of them.

As tests were perfomed on a Scientific Linux 7 system using a Intel Xeon E5-2670 CPU running at 2.60GHz. All solvers used only a single thread.

Table~\ref{allthree} shows results from 500 runs of FilterSQP
with Hessians (SQP-withHess), FilterSQP without Hessians (SQP-noHess) and our SLP
implementation (SLP-CPLEX).
Solutions are counted as feasible if the point at which the algorithm
stops has a constraint violation of $<10^{-6}$, independent of the
status returned by the solver.
Column `Ti' shows the average solution time per run, column `It' the average number of SQP or SLP iterations per run and
column `\%Good' gives the percentage of solutions that are feasible and whose
objective value is within $0.2\%$ of the best know solution from any
method (which is an acceptable tolerance in practice).

A direct comparison if SQP with either of the two SLP variants is shown in Table \ref{cfSQPSLP}.  The arrows indicate if larger or smaller numbers indicate better results. For columns `Ti', `It' and `Ti/It' we give the value for the SLP variants as a percentage of the corresponding value for SQP-withHess. Column `\%Good' gives the percentage point difference of good solutions found between the solvers (positive numbers indicating that SLP found more good solutions, negative numbers show an advantage of SQP). Generally the quality of
the solutions found by the SLP variants (SQP-noHessian and SLP-CPLEX) are better than
SQP-withHess: the
percentage of runs that are within a tolerance of 0.2\% of
the the best know is on average 4.7\% and 18.5\% higher for SQP-noHess and
SLP-CPLEX respectively.  The solution times per run are also better,
significantly so in the case of SLP-CPLEX.

The
improvements in time per run are due both to a reduction in the number
of SQP or SLP iterations and in the time per iteration.  On average
SQP-noHess and SLP-CPLEX take 55.2\% and 3.8\% of the SQP-noHess time. This
can be accounted for by the time per iteration being significantly
less (64.1\% and 11.6\%), and also because the number of iterations is
less (93\% and 35\%). SQP-withHess and SQP-noHess use the same solver,
bqpd, and the reduction in time is due to LP iterations being faster
than QP iterations. The further big reduction in time per iteration
achieved by SLP-CPLEX is due to the faster LP implementation in CPLEX
compared to bqpd.  A more surprising reason for the improved time per
run is the fact that SQP-noHess and SLP-CPLEX take fewer iterations
than SQP-withHess (average 93\% and 35\% respectively).
Intuitively SQP would be expected to be superior to SLP: after all it
uses a higher order approximation of the nonlinear programming problem
at each iteration. As one (but not the only) consequence SQP displays
second order convergence inherited from Newton's method once it has
reached a point close enough to the solution. 
SLP methods on
the other hand may have to resort to reducing the trust region radius
to zero after many steps being rejected by the filter in order to
terminate. This results in SLP often taking more iterations than SQP
to converge to high accuracy.

We do observer this tail effect in our experiments, but note that due to LP
hotstarts these iterations are very fast, often requiring only 1 or
even 0 simplex basis updates. A larger effect, however, is that away from the neighbourhood of the solution SQP is observed to
repeatedly enter the restoration phase before the algorithm is able to
“home in” on a solution, and this results in both a higher iteration
count and a decreased likelihood of finding a feasible solution
than with SLP.  


\begin{table}
\begin{center}
{\small
\begin{tabular}{l|rrr|rrr|rrr}
 & \multicolumn{3}{c|}{SQP-withHess} & \multicolumn{3}{c}{SQP-noHess} & \multicolumn{3}{c}{SLP-CPLEX} \\
problem	&	Ti	&	\%Good	&	It	&	Ti	&	\%Good	& It	&	Ti 	&	\%Good	& It 	\\
\hline
\sonsev	&	6.6	&	66.4	&	64.9	&	3.9	&	81.0	&	42.2	&	0.5	&	79.0	&	38.0	\\
\soneight	&	44.9	&	2.0	&	254.5	&	19.2	&	13.2	&	202.2	&	1.3	&	19.4	&	58.7	\\
\sept	&	20.6	&	11.4	&	105.6	&	8.0	&	6.0	&	111.6	&	1.3	&	20.0	&	50.6	\\
\jan	&	42.3	&	12.2	&	127.4	&	30.2	&	10.0	&	100.7	&	1.7	&	19.8	&	49.1	\\
\gns	&	145.5	&	20.2	&	366.1	&	26.1	&	56.0	&	340.0	&	1.7	&	99.6	&	61.2	\\
\gnsorig	&	24.2	&	97.8	&	106.2	&	13.2	&	99.2	&	151.9	&	1.5	&	100.0	&	53.9	\\
\nville	&	400.8	&	0.2	&	199.5	&	169.5	&	0.4	&	92.7	&	6.5	&	5.2	&	47.0	\\
\nvilleorig	&	112.1	&	0.0	&	110.3	&	87.7	&	0.0	&	83.1	&	4.5	&	4.0	&	47.7	\\
\mcvfor	&	681.0	&	0.4	&	359.0	&	508.0	&	0.0	&	349.0	&	14.7	&	52.6	&	48.5	\\
\mcvsevorig	&	99.0	&	32.6	&	119.1	&	58.0	&	27.2	&	161.3	&	4.6	&	20.4	&	61.5	\\
\mcvfororig	&	307.4	&	27.0	&	180.1	&	208.9	&	33.8	&	198.2	&	12.3	&	71.8	&	72.0	\\
\hline
\end{tabular}
\caption{\label{allthree}Performance of SQP-withHess, SQP-noHess and SLP-CPLEX:\\
Average per run: Ti (solution time in sec), It (number of iterations),
\%Good (percentage of solutions within 0.2\% of best known.}
}
\end{center}
\end{table}

\begin{table}
\begin{center}
{\small
\begin{tabular}{l|rrrr|rrrr}
 & \multicolumn{4}{c|}{SQP-withHess v SQP-noHess} & \multicolumn{4}{c}{SQP-withHess v SLP-CPLEX} \\
	&	Ti$\downarrow$ &	\%Good$\uparrow$ &	It$\downarrow$	&Ti/it$\downarrow$ & Ti$\downarrow$ & \%Good$\uparrow$ & It$\downarrow$	& Ti/It$\downarrow$	\\
\hline
\sonsev	&	59.1	&	14.6	&	65.0	&	90.9	&	7.58	&	12.6	&	58.55	&	12.94	\\
\soneight	&	42.8	&	11.2	&	79.4	&	53.8	&	2.90	&	17.4	&	23.06	&	12.55	\\
\sept	&	38.8	&	-5.4	&	105.7	&	36.7	&	6.31	&	8.6	&	47.92	&	13.17	\\
\jan	&	71.4	&	-2.2	&	79.0	&	90.3	&	4.02	&	7.6	&	38.54	&	10.43	\\
\gns	&	17.9	&	35.8	&	92.9	&	19.3	&	1.17	&	79.4	&	16.72	&	6.99	\\
\gnsorig	&	54.5	&	1.4	&	143.0	&	38.1	&	6.20	&	2.2	&	50.75	&	12.21	\\
\nville	&	42.3	&	0.2	&	46.5	&	91.0	&	1.62	&	5.0	&	23.56	&	6.88	\\
\nvilleorig	&	78.2	&	0.0	&	75.3	&	103.8	&	4.01	&	4.0	&	43.25	&	9.28	\\
\mcvfor	&	74.6	&	-0.4	&	97.2	&	76.7	&	2.16	&	52.2	&	13.51	&	15.98	\\
\mcvsevorig	&	58.6	&	-5.4	&	135.4	&	43.3	&	4.65	&	-12.2	&	51.64	&	9.00	\\
\mcvfororig	&	68.0	&	6.8	&	110.0	&	61.8	&	4.00	&	44.8	&	39.98	&	10.01	\\
\hline
Average	&	55.2	&	4.7	&	93.1	&	64.1	&	3.81	&	18.5	&	35.21	&	10.56	\\
\hline
\end{tabular}
\caption{\label{cfSQPSLP}Improvements from SQP-withHess to SQP-noHess or SLP-CPLEX:\\
Ti$\downarrow$, It$\downarrow$ and Ti/It$\downarrow$ are the ratio (as \%) of SQP-withHess values to SQP-noHess or SLP-CPLEX values. \%Good$\uparrow$ is the difference between the  SLP-CPLEX or SQP-noHess quality and the SQP-withHess quality.}
}

\end{center}
\end{table}

As an explanation of this behaviour we offer the following insight:
the nonlinearity in the pooling problem is exclusively due to
bilinear terms; these give indefinite a Hessian contribution of the
form
$$
\left[\begin{array}{cc} 0&1\\1&0
\end{array}\right].
$$
In fact rather than being helpful, these Hessians bias the algorithm
towards taking steps along negative curvature directions (increase one
bilinear variable while decreasing the other, as much as possible)
which is not desirable. Away from the region of quadratic convergence
of Newtons method (that is, for the vast majority of the SQP
iterations) this can lead to rather erratic behaviour of SQP.

When comparing the two SLP variants, SQP-noHess and SLP-CPLEX, it can be observed that SLP-CPLEX takes significantly fewer
iterations and less time to solve each problem, but most noticeably
massively increases the likelihood of finding a feasible solution.
This does not seem well explained by the algorithmic differences
between them. In fact, while the main reason for
an infeasible run in the plain SQP algorithm, is convergence to a local
solution of the feasibility (phase-I) problem, in the no-Hessian
version the main reason for infeasibility is algorithmic failure due
to inconsistent second order information -- often there are long
sequences of rejected second order correction steps that reduce the
trust region radius and subsequently lead to premature termination of
the algorithm. In SLP-CPLEX this inconsistent algorithm logic has been
removed.

\begin{table}
\begin{center}
{\small
\begin{tabular}{l|rrr|rrr|rrr}
 & \multicolumn{3}{c|}{SLP-CPLEX} & \multicolumn{3}{c}{SQP-withHess} & \multicolumn{3}{c}{IPOpt} \\
problem	&	Ti	& \%Feas & \%Good & Ti & \%Feas & \%Good & Ti & \%Feas & \%Good \\
\hline
\sonsev	&	0.5	&	100.0	&	79.0	&	6.6	&	80.6	&	66.4	&	6.4	&	96.2	&	93.2	\\
\soneight	&	1.3	&	100.0	&	19.4	&	44.9	&	12.6	&	2.0	&	30.9	&	0.0	&	0.0	\\
\sept	&	1.3	&	99.2	&	20.0	&	20.6	&	96.2	&	11.4	&	33.9	&	50.4	&	33.8	\\
\jan	&	1.7	&	98.5	&	19.8	&	42.3	&	56.4	&	12.2	&	46.1	&	53.4	&	12.8	\\
\gns	&	1.7	&	100.0	&	99.6	&	145.5	&	20.6	&	20.2	&	58.8	&	85.0	&	86.4	\\
\gnsorig	&	1.5	&	99.8	&	100.0	&	24.2	&	98.8	&	97.8	&	67.3	&	33.2	&	33.2	\\
\nville	&	6.5	&	95.0	&	5.2	&	400.8	&	15.4	&	0.2	&	71.3	&	99.4	&	9.6	\\
\nvilleorig	&	4.5	&	78.4	&	4.0	&	112.1	&	42.8	&	0.0	&	37.6	&	99.4	&	6.8	\\
\mcvfor	&	14.7	&	99.0	&	52.6	&	681.0	&	22.4	&	0.4	&	508.6	&	58.2	&	32.0	\\
\mcvsevorig	&	4.6	&	99.0	&	20.4	&	99.0	&	81.0	&	32.6	&	82.6	&	0.0	&	0.0	\\
\mcvfororig	&	12.3	&	98.8	&	71.8	&	307.4	&	72.6	&	27.0	&	437.7	&	2.2	&	0.6	\\
\hline
\end{tabular}
\caption{\label{tab:slpsqp}Performance of SLP-CPLEX, SQP-withHess and IPOpt:\\
Ti (average time per run in sec),
\%Feas  (percentage of runs that are feasible),\\
\%Good (percentage of runs with a solutions within 0.2\% of best known)
.}
}
\end{center}
\end{table}

Table~\ref{tab:slpsqp} compares the average solution time of SLP-CPLEX, SQP (with Hessians) and IPOpt. Generally IPOpt takes about the same amount of time as SQP, although with some large variation. In all cases, however, SLP is an order of magnitude faster than either of the other two algorithms. Somewhat surprisingly IPOpt struggles to find feasible solutions for some problems, in particular \soneight, \mcvsevorig~ and \mcvfororig~where none or almost none of the runs where feasible. 

\begin{table}
\begin{center}
{\small
\begin{tabular}{l|rrrr|rrr}
& \multicolumn{4}{c|}{ETiGood} & \multicolumn{3}{c}{ SLP-CPLEX Speedup relative to:} \\
Problem & SLP-CPLEX & SQP-withHess & IPOpt & Baron & SQP-withHess & IPOpt & Baron \\
\hline
\sonsev	&	0.6	&	9.9	&	6.9	&	3	&	15.7	&		10.8	&	4.7	\\
\soneight	&	6.7	&	2245.0	&	inf	&	11	&	335.0	&		inf	&	1.6	\\
\sept	&	6.5	&	180.7	&	100.3	&	20	&	27.7	&		15.4	&	3.1	\\
\jan	&	8.6	&	346.7	&	360.2	&	21	&	40.4	&		41.9	&	2.4	\\
\gns	&	1.7	&	721.7	&	68.1	&	151	&	422.8	&		39.9	&	88.5	\\
\gnsorig	&	1.5	&	24.7	&	202.7	&	51	&	16.5	&		135.1	&	34.0	\\
\nville	&	125.0	&	200400.0	&	742.7	&	2331	&	1603.2	&		5.9	&	18.6	\\
\nvilleorig	&	112.5	&	inf	&	552.9	&	147	&	inf	&		4.9	&	1.3	\\
\mcvfor	&	27.9	&	170250.0	&	1589.4	&	1362	&	6091.9	&		56.9	&	48.7	\\
\mcvsevorig	&	22.5	&	303.7	&	inf	&	268	&	13.5	&	inf	&	11.9	\\
\mcvfororig	&	17.1	&	1138.5	&	72950.0	&	107	&	66.5	&		4258.4	&	6.2	\\
\hline
Average	&		&		&		&		&	863.3	&		507.7	&	20.1	\\
\hline
\end{tabular}
\caption{\label{cfIPOPTSQPSLP}
ETiGood is the expected time to find a Good ({\em i.e.} 0.2\%) solution.\\
EIiGood = $\mbox{Ti}{100\over \mbox{\%Good}}$ 
.}
}
\end{center}
\end{table}

Another way to look at the results in Table 4 would be in terms of
 {\em Expected Time to find a Good solution} (ETiGood) which can be worked out as  the
ratio 100Ti/\%Good. These are presented in Table~\ref{cfIPOPTSQPSLP}.
On this measure
SLP-CPLEX is significantly faster than the others.
The ETiGood speed up of SLP-CPLEX relative to SQP-withHess is in the range
13.5 to 6,091.9 with average 863.3 (excluding the 1 problem where SQP-withHess
failed to find a Good solution).
The ETiGood speed up of SLP-CPLEX relative to IPOPT is in the range
4.9 to 4,258.4 with average 570.6 (excluding the 2 problem where IPOpt
failed to find a Good solution).
The ETiGood speed up of SLP-CPLEX relative to Baron is in the range
1.3 to 88.5 with average 20.1.

\subsection{Comparison of local and global solvers}

We have compared the performance of the global solver Baron using the
formulations presented in Section~\ref{sec:form} with the two best local solvers, namely SLP-CPLEX and IPOpt.  
Each solve employed a 2 hour time limit
for Baron. For the local solvers we have only used the qq-formulation. The results of these experiments are presented in
Figure~\ref{fig:s97}.

The graphs should be read as follows: Each
gives a plot of quality of solution found ($y$-axis) vs time spend in seconds
($x$-axis, logarithmic). For Baron (for each of the four different formulations)
this is a straightforward plot of the progression of the best feasible
solution found within the time limit. The lower bounds on the solution obtained are not shown as they are much lower and off the scale of most of the graphs.
Indeed the {\tt qq+}, {\tt pq} and {\tt pqs} formulations all
  obtained the same lower bound (at the root node) that could not be
  improved within the 2h time limit 
The
  qq-formulation leads to significantly weaker lower bounds (gap larger by a
  factor of 10). The lower bounds are not shown on the graphs since they are (sometimes significantly) below the lower limit of the area shown. Values can be seen in Table~\ref{tab:testproblems}.

For SLP (blue curve) and IPOpt (red curve) on the qq-formulation we show the expected time
that would be needed to obtain a solution at least as good as a given
objective value $\hat{v}$. That is let $\bar{t}$ be the average time
needed for a local solve.
If $k$ out of $n$ SLP/IPOpt runs find a solution better than $\hat{v}$ we model this as a Bernoulli trial with success probability $p = k/n$. The expected number of trials until the first success is $1/p = n/k$ runs or time $\hat{t} = \bar{t}n/k$. 
From the figures it can be seen that 
\begin{enumerate}
\item The SLP and IPOpt curves are almost smooth (rather than step functions) indicating that the problems have a huge number of local optima. This ties in with the earlier analysis in Section 3.
\item Baron always finds better solutions with the qq-formulation than with the
pq-formulation. The main reason seems to be that the qq-formulation is
smaller (since it does not include the pq-constraints) and thus is able to process many more nodes in the same
time. In fact the first feasible solution is found by the
qq-formulation much faster than for the pq-formulation.
\item Strengthening the qq-formulation by the pq-cut ({\tt qq+}) does not pay off: neither for the local solvers (due to the larger problem size), nor, somewhat surprisingly, for Baron.
  While it does strengthen the lower bound,
  the resulting increase in problem size means that it takes 
  longer to find solutions of the same quality.
\item However, with Baron, the qq-formulation strengthened by the pq-cut ({\tt qq+}) performs better than the pq-formulation (to which it is in some sense equivalent). 
\item[5a.] SLP is clearly superior to Baron in terms of time taken to find a solution of a given quality: The SLP curve lies well below the Baron curve for almost all problems, times and formulations.
\item[5b.] Only for the qq-formulation there are a few problems ({\tt \sonsev, \soneight, \nvilleorig, \gns, \gnsorig}) where the first feasible solution found by Baron is better than the best solution that could
  be expected to be found by SLP in the same time 
and only for {\tt \soneight} is the difference more than marginal. However given more
  time SLP will find a better solution. Also
  SLP will have already have found very good solutions before Baron with the
  qq-formulations has found the first feasible solution.
\item[6.] At the end of the 2h time limit SLP has always found the best
  known solution (as given in Table~\ref{tab:testproblems}). 
Baron with the qq-formulation fails to find the best solution within 2h for problems {\tt \sonsev, \gns, \gnsorig, \mcvfor} and {\tt \mcvfororig}. 
\item[7.] The curve of IPOpt looks similar to the one for SLP but shifted to the right. The reason is that IPOpt needs much longer for a single run ($\overline{t}$) than SLP does. The IPOpt curves seem to drop down quicker than the SLP ones. This is mainly due to the logarithmic scale of the $x$-axis, but also for some problems due to the probability of finding a solution close to globally optimal being larger for IPOpt than for SLP-CPLEX (see Table~\ref{tab:slpsqp}). IPOpt is however clearly uncompetitive for all problems. For problems {\tt \soneight, \mcvsevorig} and {\tt \mcvfororig} the IPOpt curve is not shown (or off the plot) since all (or almost all) of the runs are infeasible.

\end{enumerate}
\section{Conclusions\label{sec:conclusions}}
We have given a comparison of several local solvers employing 
randomized starting points with that of the global solver Baron for the
pooling problem. 

The best local solver is sequential linear programming, which -- while the
oldest method -- somewhat surprisingly significantly outperforms newer
methods such as SQP and Interior Point. 
Measured by best quality solution obtained in a given
time we find that SLP performs much better than Baron for all problem formulations, including the traditional pq-formulation, and almost all time limits. 
In terms of the expected time to find a solution within 0.2\% of the best known SLP with CPLEX as subproblem solver shows an average speedup of respectively 863, 507 and 20 times relative to FilterSQP, IPOpt and Baron.

We further propose a new formulation of the pooling problem, which we
term the qq-formulation. This is of comparable size to the q-formulation and can be strengthened by additional cuts analoguous to pq-cuts. 
When measured as quality of solution found in
a given time Baron's performance on the qq-formulation is superior to all other formulations. However, strengthening this
formulation with the pq-like-cuts is not worthwhile on this criterion. 
The qq-formulation is the only one for which Baron is not totally dominated
by SLP for all time limits. With the qq-formulation for a minority of problems 
the first feasible solution found by Baron is
marginally better than the solution that could be expected by randomized SLP in
the same time. However for the vast majority of time limits SLP returns better solutions than Baron even for the qq-formulation.
One advantage of Baron is that it is able to return a lower bound (and
thus an optimality gap), which is not the case for any of the local
solvers. However for most of the test examples the bound gap achieved
is too large to be of practical value to the
problem owner.



\begin{figure}[p!]
\centering
  \begin{tabular}[b]{@{}c@{}c@{}}
   \HGAP{-0.05}\includegraphics[width=7.5cm]{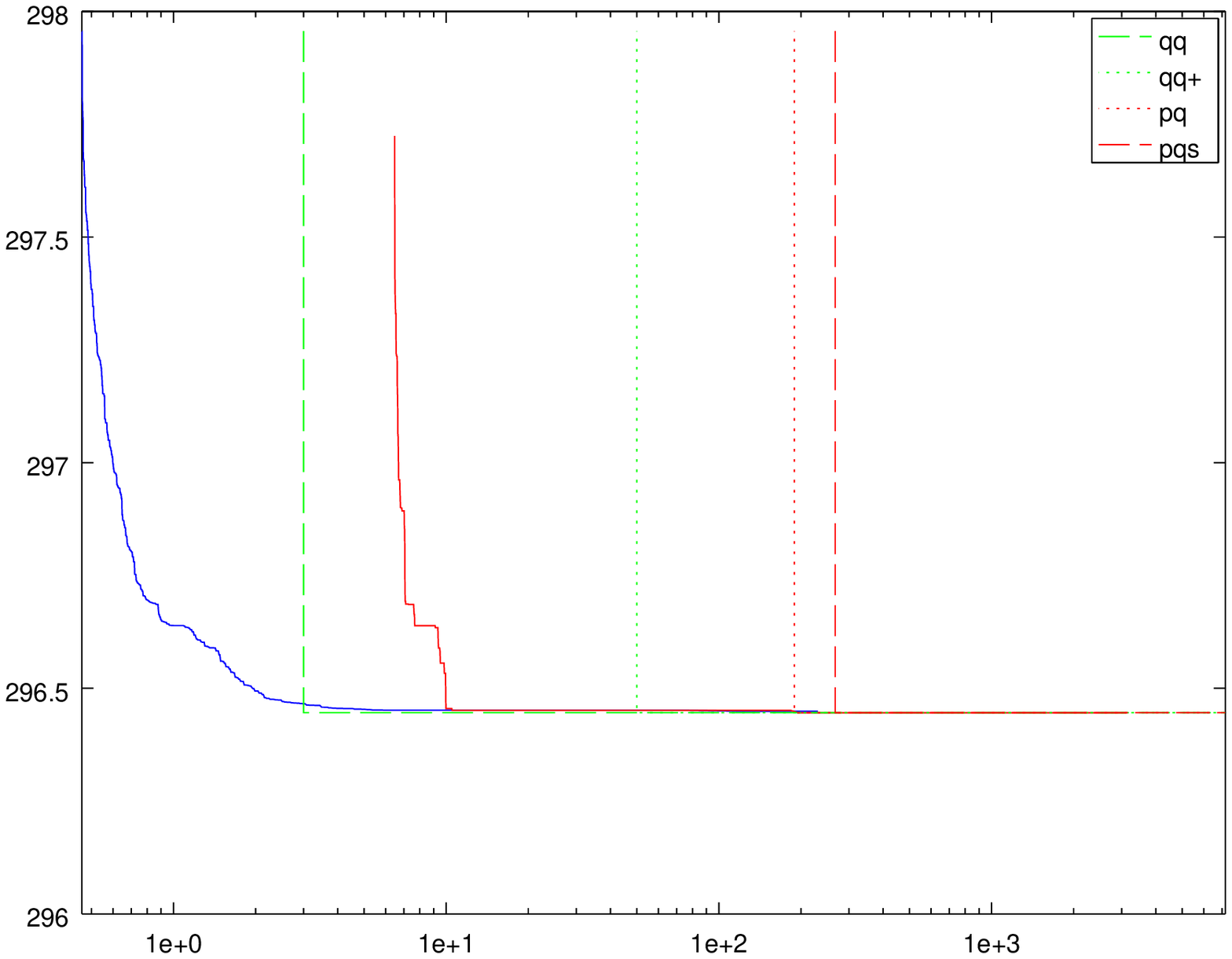} &
    \HGAP{-0.02}\includegraphics[width=7.5cm]{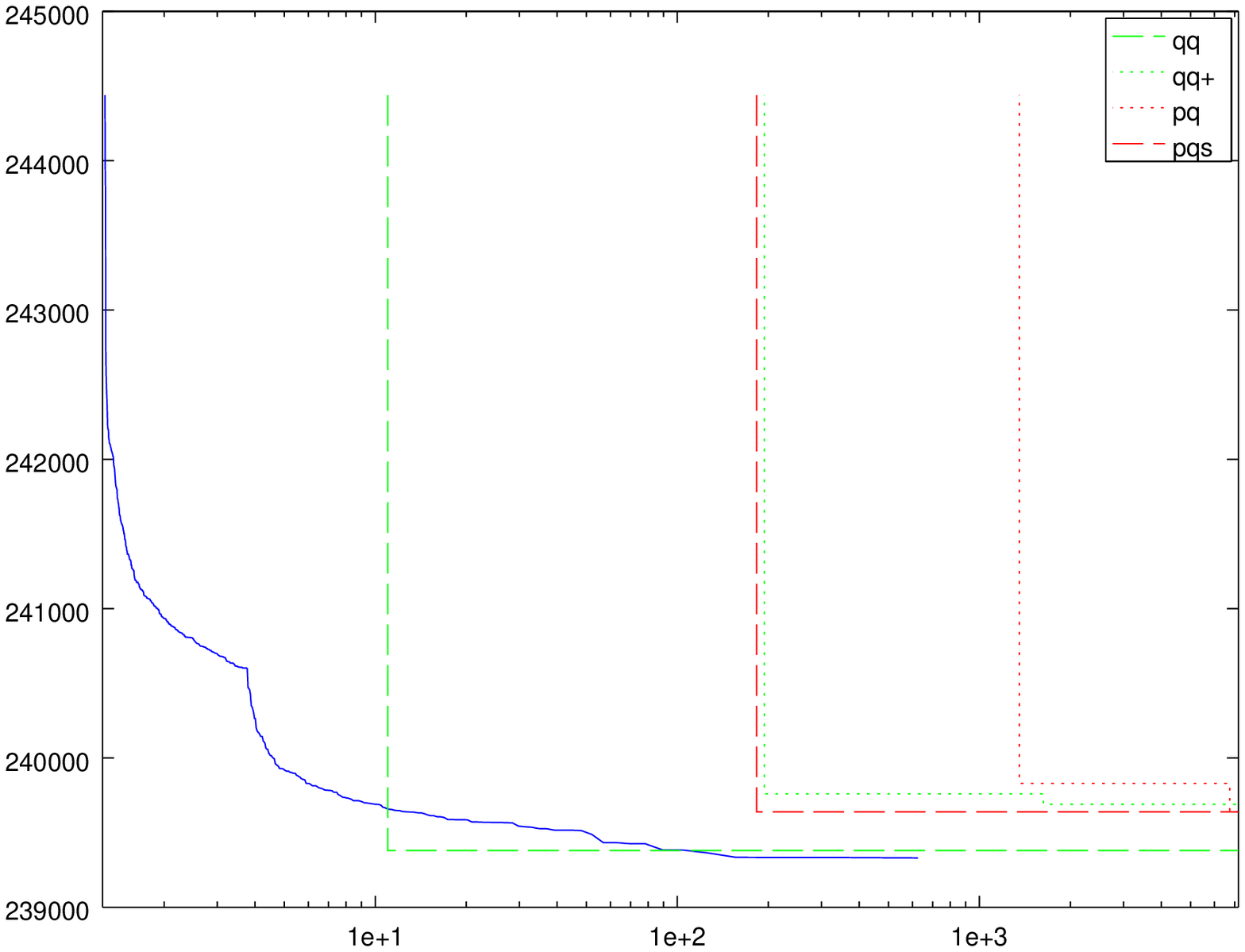} \\
    \small (a)\sonsev&   \small (b)\soneight
  \end{tabular}

\centering
  \begin{tabular}[b]{@{}c@{}c@{}}
    \HGAP{-0.05}\includegraphics[width=7.5cm]{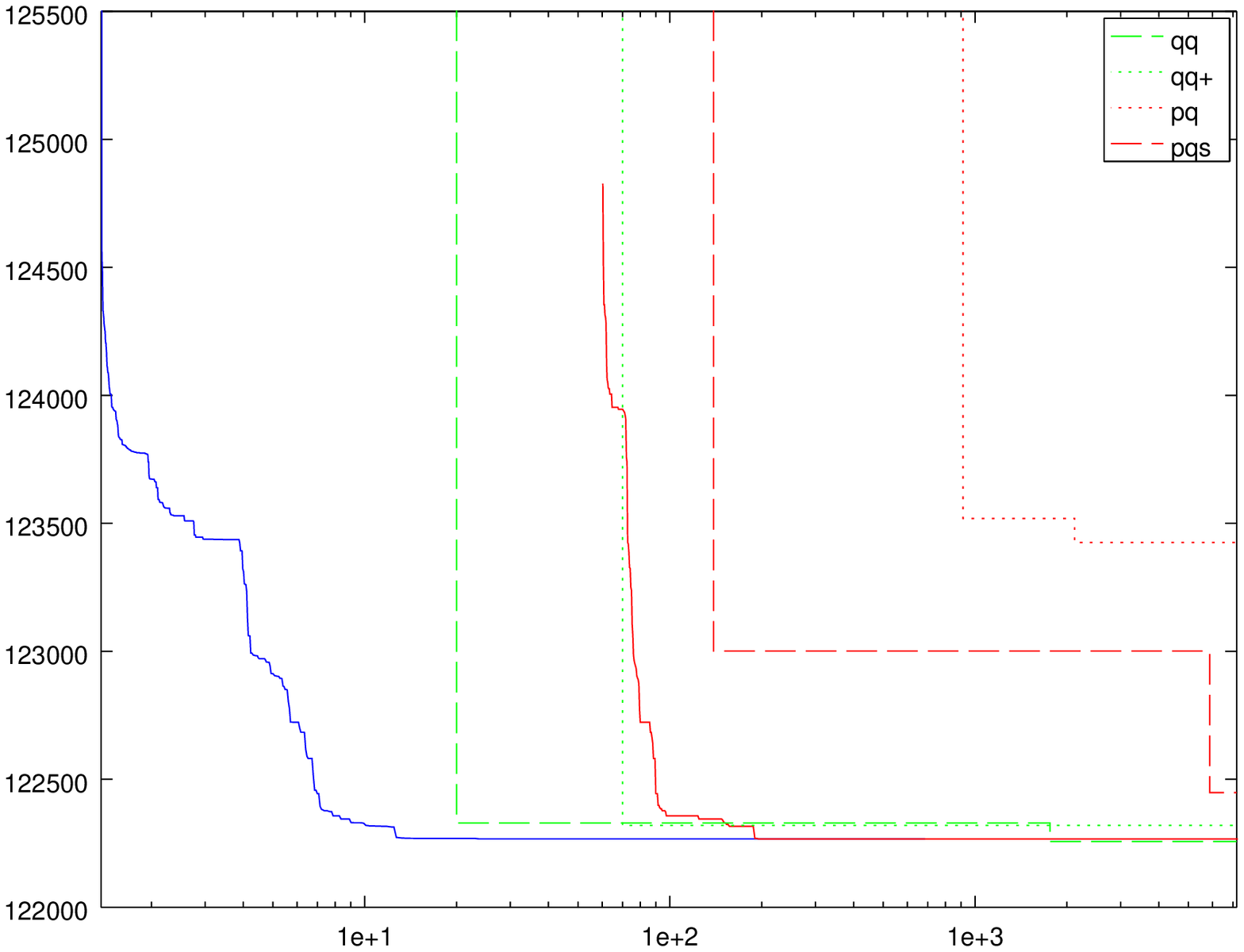} &
    \HGAP{-0.02}\includegraphics[width=7.5cm]{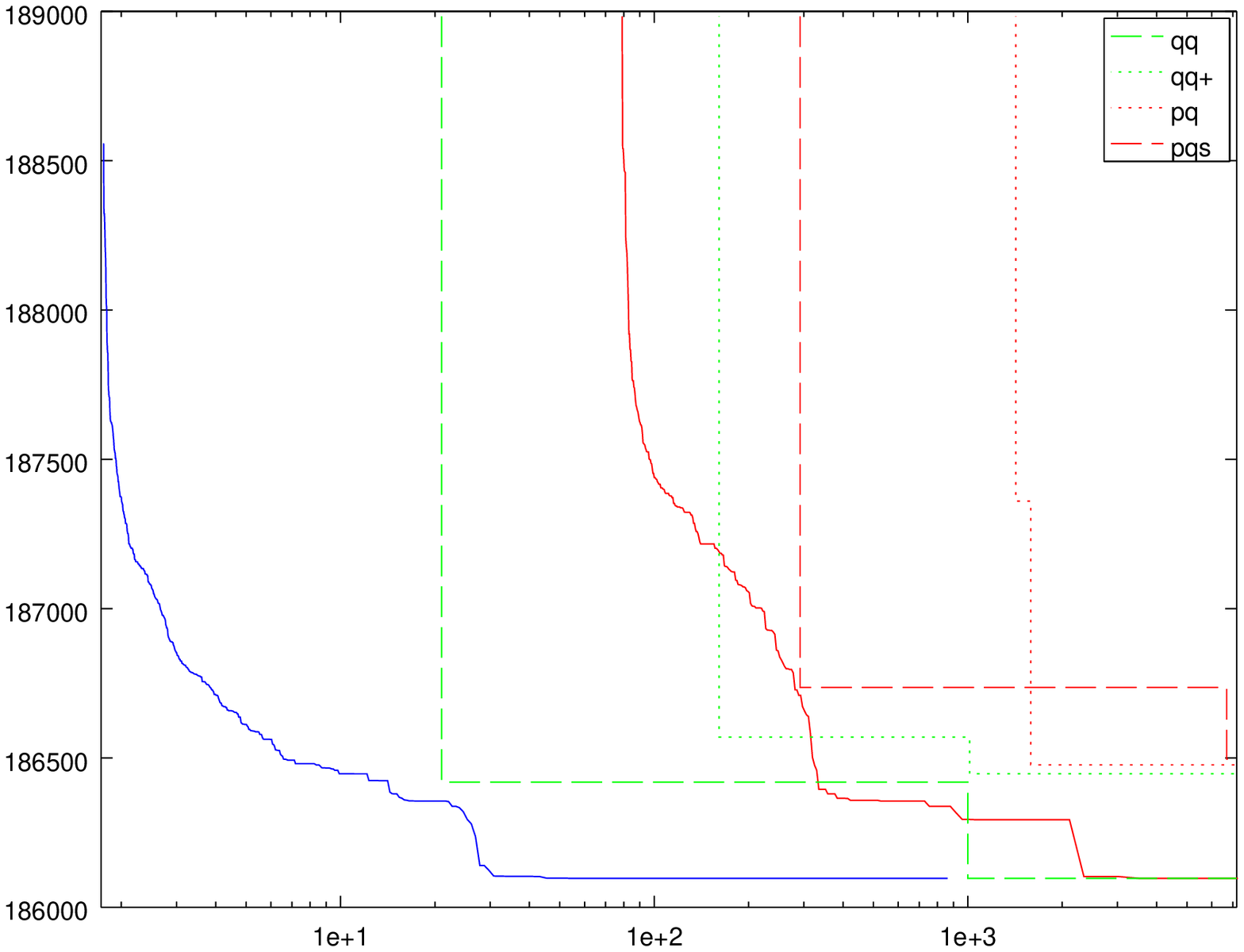} \\
    \small (c)\sept&   \small (d)\jan
  \end{tabular}

\centering
  \begin{tabular}[b]{@{}c@{}c@{}}
    \HGAP{-0.05}\includegraphics[width=7.5cm]{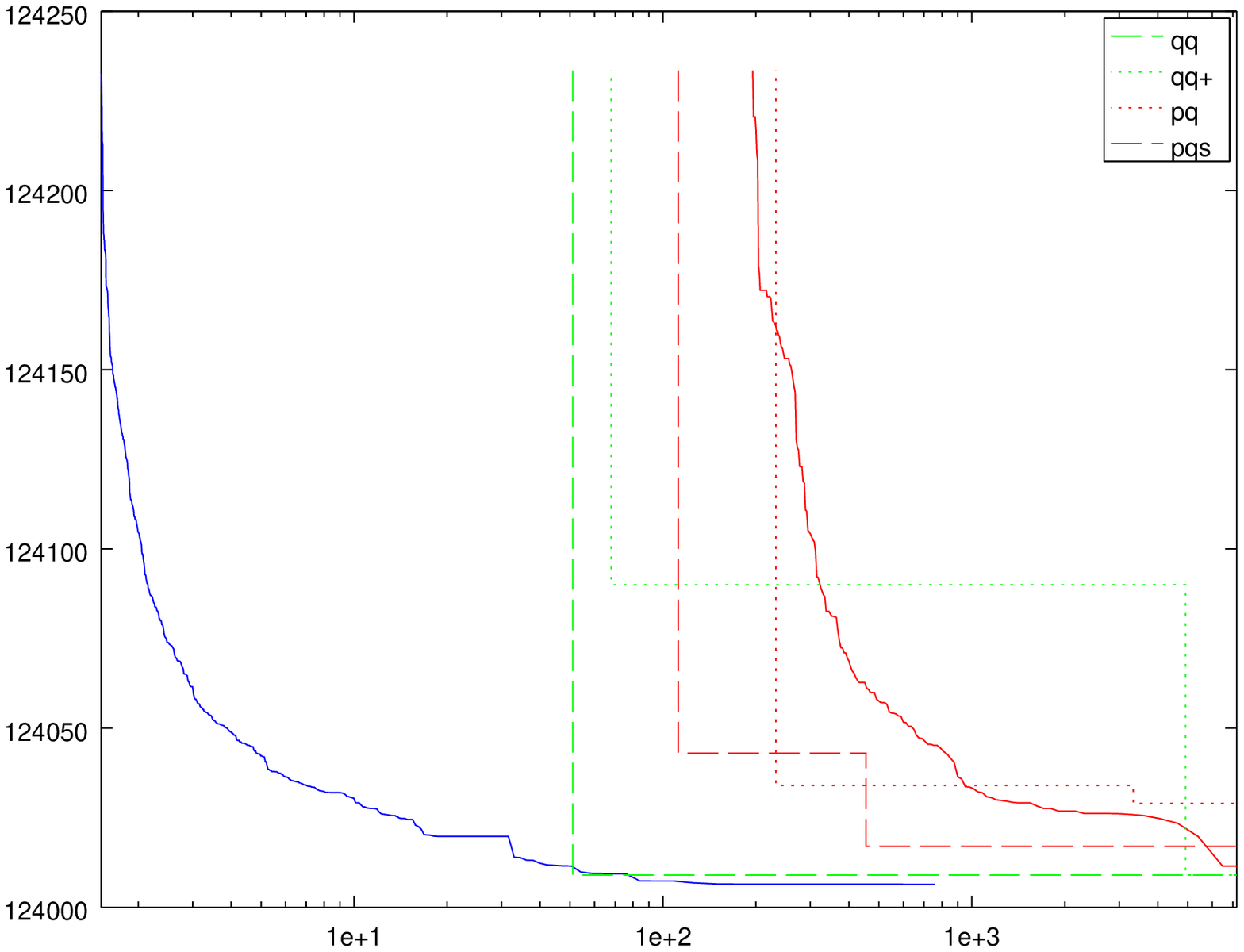} &
    \HGAP{-0.02}\includegraphics[width=7.5cm]{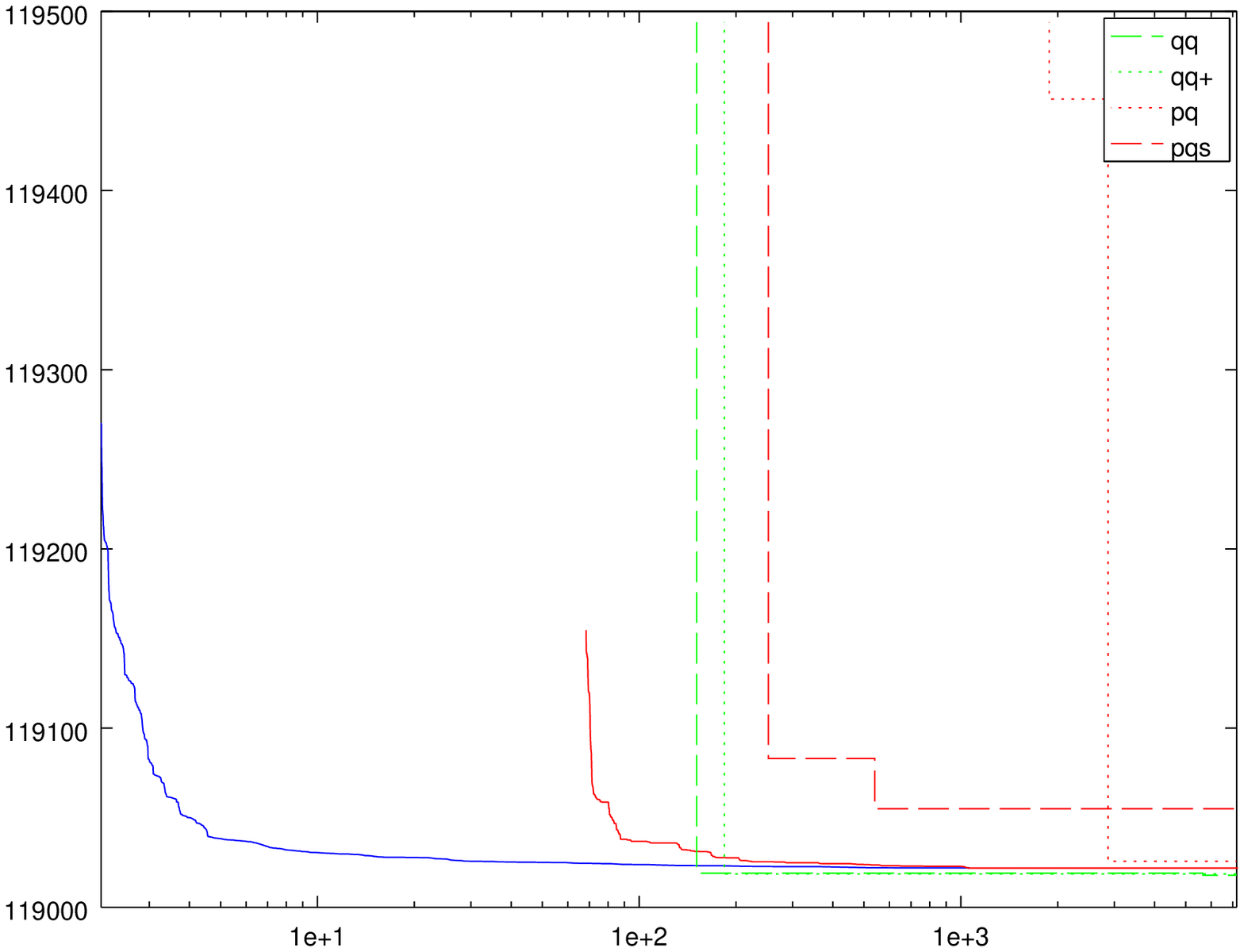} \\
    \small (e)\gnsorig&   \small (f)\gns
  \end{tabular}
\caption{\label{fig:s97}Comparing SLP (blue) with IPOpt (red) and Baron for different formulations.}
\end{figure}

\begin{figure}[p!]
\setcounter{figure}{3}
\centering
  \begin{tabular}[b]{@{}c@{}c@{}}
    \HGAP{-0.05}\includegraphics[width=7.5cm]{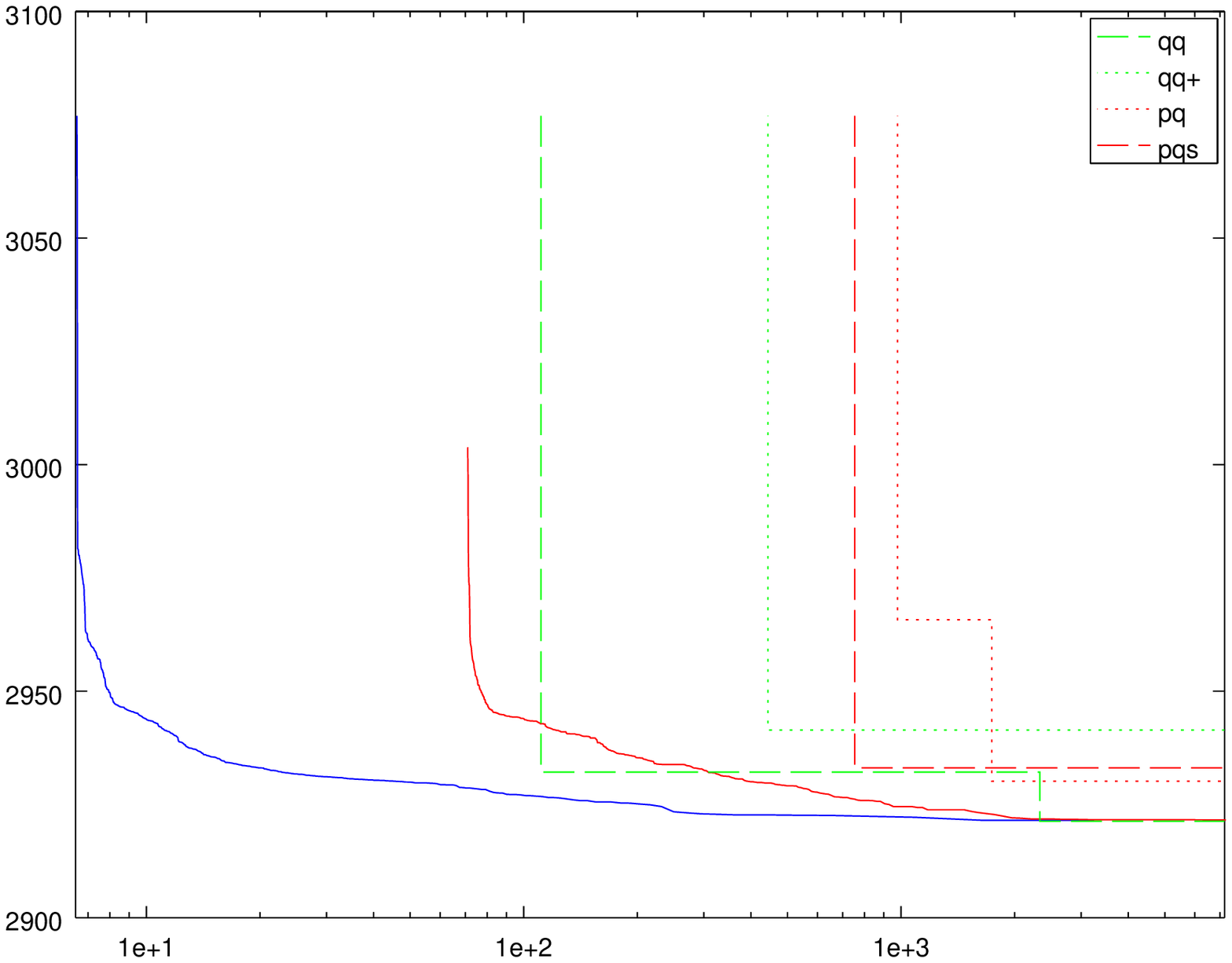} &
    \HGAP{-0.02}\includegraphics[width=7.5cm]{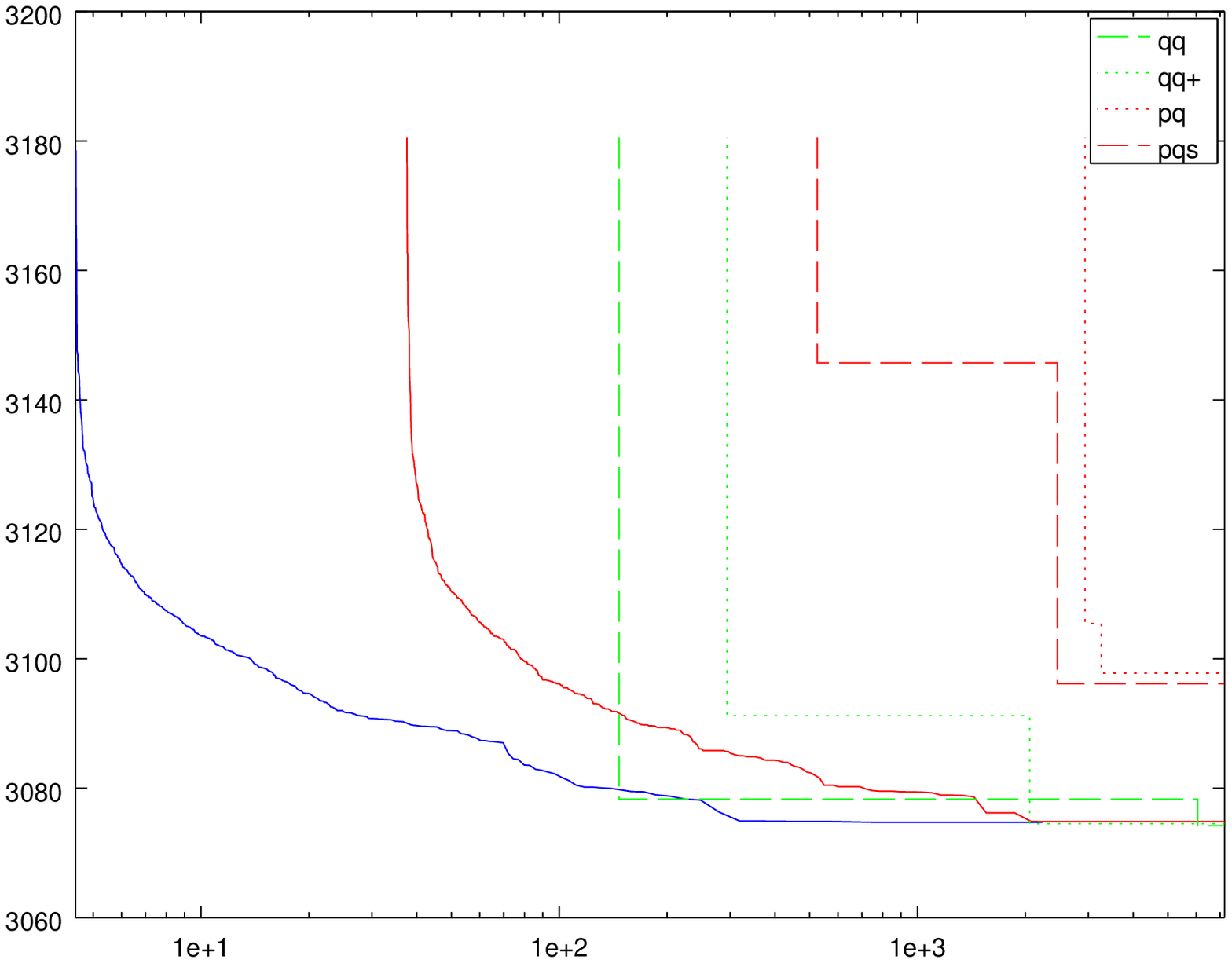} \\
    \small (g)\nville&   \small (h)\nvilleorig
  \end{tabular}

%
\centering
  \begin{tabular}[b]{@{}c@{}c@{}}
    \HGAP{-0.05}\includegraphics[width=7.5cm]{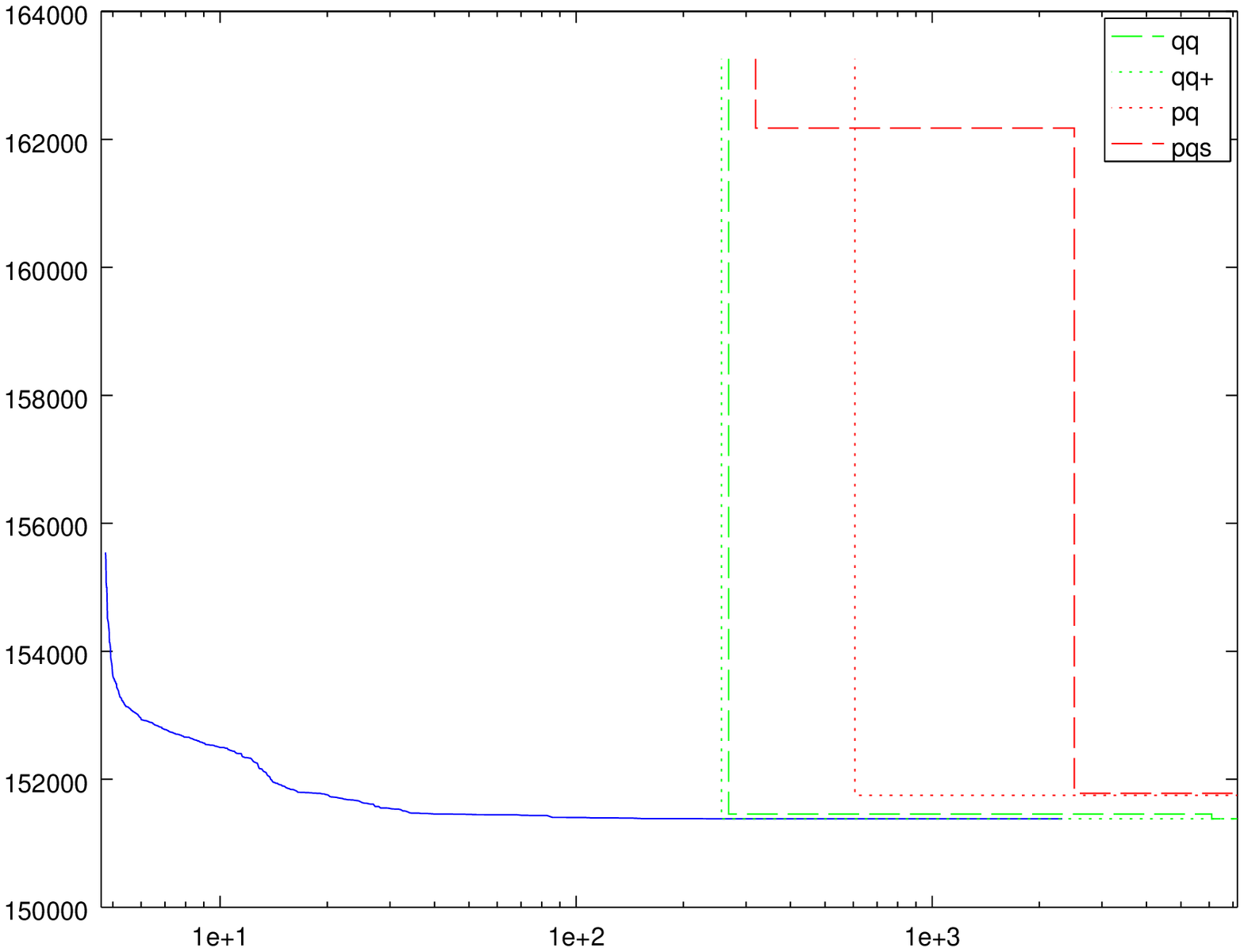} &
    \HGAP{-0.02}\includegraphics[width=7.5cm]{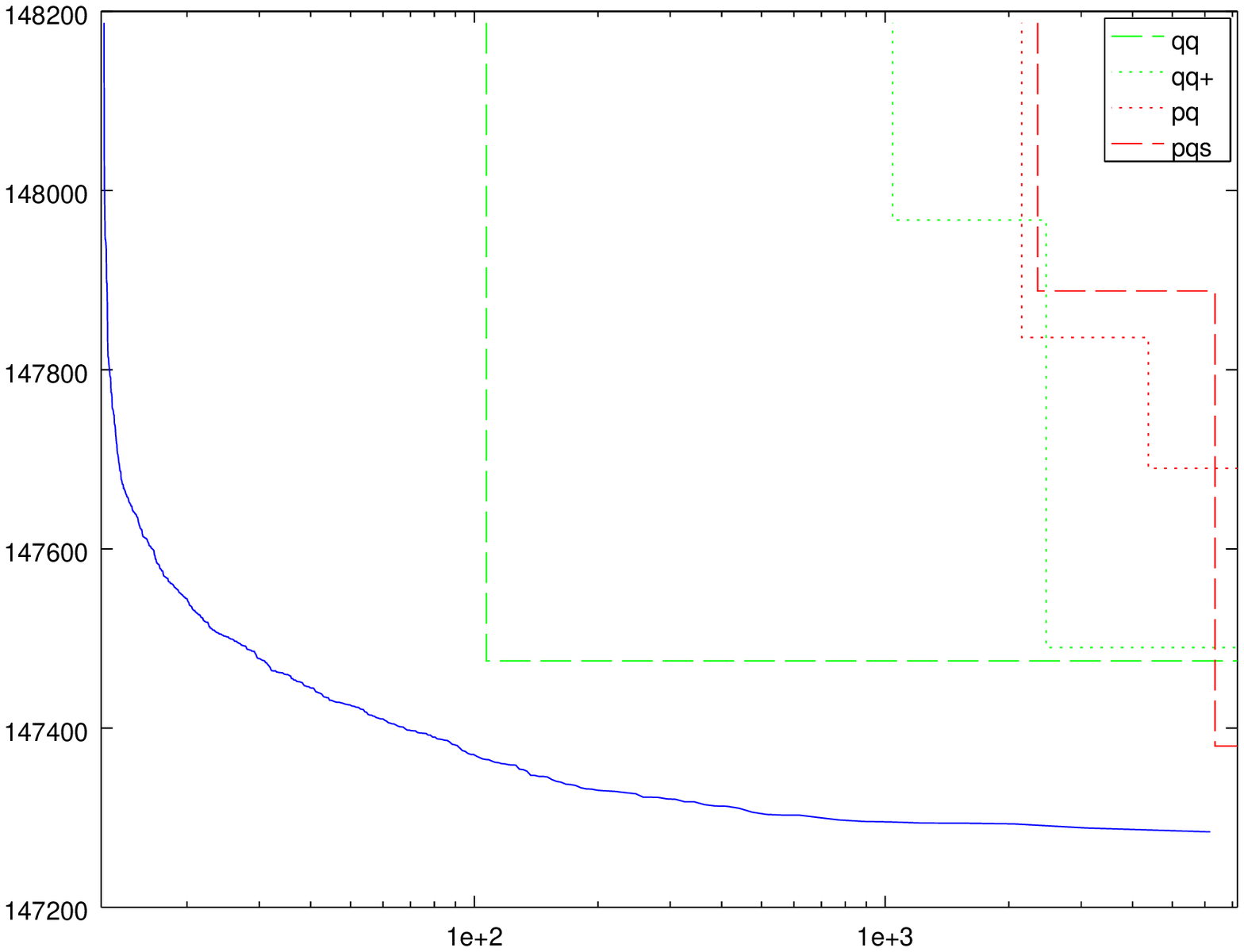} \\
    \small (i)\mcvsevorig&   \small (j)\mcvfororig
  \end{tabular}

\centering
  \begin{tabular}[b]{@{}c@{}c@{}}
    \HGAP{-0.05}\includegraphics[width=7.5cm]{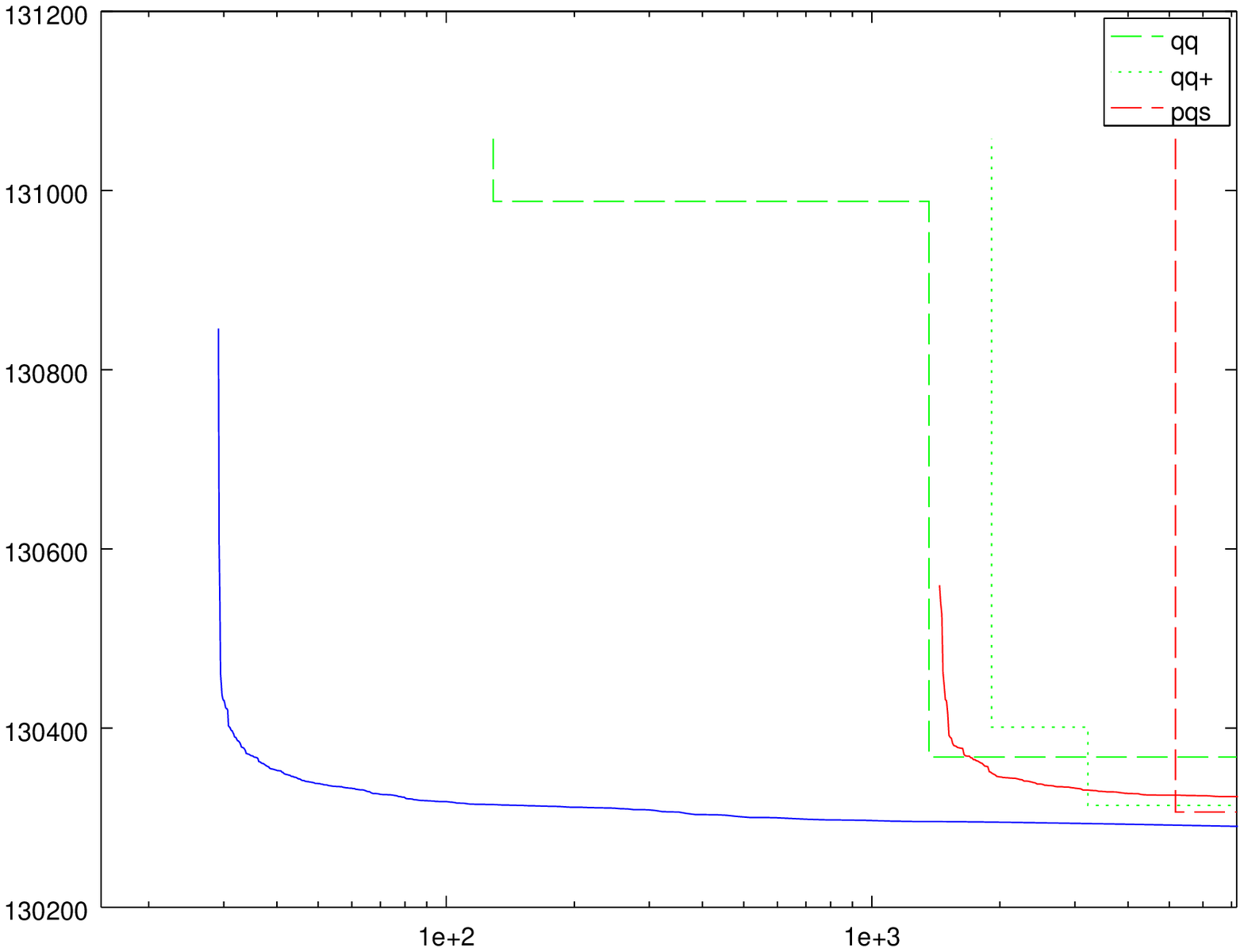} \\
    \small (k)\mcvfor
  \end{tabular}
\caption{Comparing SLP (blue solid) with IPOpt (red solid) and Baron for different formulations (cont).}
\end{figure}

\bibliographystyle{plain}
\bibliography{pooling}

\end{document}